\documentclass[preprint, 3pt]{elsarticle}

\RequirePackage{amsmath}

\usepackage[english]{babel}
\usepackage[T1]{fontenc}
\usepackage{empheq}
\usepackage{amsfonts}
\usepackage{enumerate}
\usepackage{amsmath}
\usepackage{amsthm}
\usepackage{mathtools}
\allowdisplaybreaks
\mathtoolsset{showonlyrefs}			
\numberwithin{equation}{section}

\newtheorem{thm}{Theorem}[section]
\newtheorem{prop}[thm]{Proposition}

\newtheorem{cor}[thm]{Corollary}
\newtheorem{lem}[thm]{Lemma}
\newtheorem{rem}[thm]{Remark}
\newtheorem{hyp}[thm]{Hypothesis}

\newcommand{\probSpace}{(\Omega,\mathcal{F},\{\mathcal{F}_{t}\}_{t\geq 0},\mathbb{P})}

\def\R{\mathbb R}
\def\N{\mathbb N}

\def\E{\mathbb E}
\def\P{\mathbb P}
\def\Q{\mathbb Q}

\def\shd{{\cal D}}
\def\shf{{\cal F}}
\def\shg{{\cal G}}

\def\shi{{\cal I}}

\def\shn{{\cal N}}
\def\shl{{\cal L}}

\def\shq{{\cal Q}}

\def\shy{{\cal Y}}

\def\shx{{\cal X}}

\def\1{\mathds{1}}
\def\ra{\rightarrow}

\def\ve{\varepsilon}
\def\a{\alpha}
\def\s{\sigma}

\def\pim{\frac{\pi}{2}}

\def\w{\wedge}

\def\um{\frac{1}{2}}

\usepackage{geometry}
\begin{document}

\begin{frontmatter}

\title{Pointwise explicit bound for derivative of solutions to linear parabolic PDEs with Neumann boundary conditions}

\author[1]{Carlo Ciccarella}
\ead{carlo.ciccarella@hotmail.com}

\affiliation[1]{organization={Institut de Mathématiques, Ecole Polytechnique Fédérale de Lausanne},
addressline={Station 8},
postcode={CH-1015},
city={Lausanne},
country={Switzerland}}

\begin{abstract}

We derive explicit pointwise bounds for the spatial derivative $\left| \frac{\partial V}{\partial x} \right|$ of solutions to linear parabolic PDEs with Neumann boundary conditions.
The bound is  fully explicit in the sense that it depends only on the coefficients of the PDE and the domain, including closed-form expression for all constants. 
The proof is purely probabilistic.  We first extend to time inhomogeneous diffusions a result concerning the derivative of the solution of a reflected SDE. Then, we combine it with the spectral expansion of the law of the first hitting time to a boundary for a reflected diffusion.
The main motivation comes from optimal control where, in order to apply verification theorems,  precise gradient estimates are often required when closed-form solutions of the Hamilton–Jacobi–Bellman equation are unavailable.
This result will be used in a forthcoming work to rigorously prove that the conjectured optimal strategy for the sailboat trajectory optimization problem  is indeed optimal far from the buoy.  
We also state a sufficient condition for $\limsup_{t\ra \infty} \left| \frac{\partial V}{\partial x}(t,x) \right|$ to be bounded,  which only involves the coefficients of the problem and the first eigenvalue of the spectral expansion. 
\end{abstract}


\end{frontmatter}

\section{Introduction}


Let $\shi$ be an open bounded interval of the  real line. For simplicity, we set $\shi = ]0,\zeta[$ for some $\zeta>0$. We denote by $\bar \shi = [0,\zeta]$ its closure,  and by $\partial \shi = \{0,\zeta \}$ its boundary. In this work, we will consider a classical solution, in a sense specified by Proposition \eqref{existSol}, to the following parabolic PDE with Neumann boundary conditions:
\begin{align}
 \frac{\partial V(t,x)}{\partial t} &=  b(t,x) \frac{\partial V(t,x)}{\partial x}  + \frac{\sigma^2}{2 } a^2(x)  \frac{\partial^2 V(t,x)}{\partial x^2}  + h(t,x)  
 \label{freeB2aa}\\
\frac{\partial V(t,0)}{\partial x} &= n(t), \qquad \frac{\partial V(t,\zeta)}{\partial x} = m(t), \,  &t \geq 0 , \label{freeB2cc} \\
V(0,x) &= g(x), \,&x \in \bar \shi \label{freeB2ee}.
\end{align}
More details on the coefficients are specified in Hypothesis~\ref{hyp}. 

\begin{rem}\label{chg_var}
By introducing the  change of variables  $U(t,x) = V(t,x) + \phi(t,x)$, where $\phi(t,x) = -n(t)x -\frac{m(t)-n(t)}{2\zeta}x^2$, the problem  \eqref{freeB2aa}--\eqref{freeB2ee} can be reformulated with homogeneous boundary conditions. In other words, the boundary terms reduce to $n=m \equiv 0$. Hence, without loss of generality, all results will be stated and proved under this latter assumption.
\end{rem}

This work, to the best of our knowledge, is the first in which a pointwise bound in closed form on the $x$-derivative of the solution to a parabolic PDE with non-constant coefficients is provided.
It is the first of a series of works, in which we will also treat the Dirichlet problem and extensions to multi-dimensional settings.

For simplicity of exposition, we assume that the PDE has no potential term, that is in \eqref{freeB2aa} there is no component of the form $f(t,x)V(t,x)$. However, the results we establish are readily generalizable to the case where such a potential is present.

We now outline the main ideas of the paper and state the principal result.

Let $(B_t)_{t\geq 0}$, a standard real-valued Brownian motion defined on a filtered probability space $\probSpace$ satisfying the usual assumptions and consider, for $x\in \bar \shi$ and $t>0$, the reflected diffusion $X^{t,x}_s$: 
\begin{align}\label{FC}
&dX_s^{t,x} = b(t-s, X_s^{t,x}) ds +\sigma \, a(X_s^{t,x})  dB_s + dL^{t,x}_s - dU^{t,x}_s, \, s\in [0, t]\\
&X_s^{t,x} \in \bar \shi , \qquad  X^{t,x}_0 = x,\\
&1_{\{X^{t,x}_s \neq 0 \}}dL^{t,x}_s=0,\qquad 1_{\{X^{t,x}_s \neq \zeta \}}dU^{t,x}_s=0.
\end{align}
Existence and uniqueness  of a  strong solution to \eqref{FC} follows  from \cite{tanaka1979}, Theorem 4.1.

The proof of our main result, Theorem~\ref{final_thm}, relies on the following key arguments.

\begin{enumerate}[(1)]
 \item 
We prove the Feynman-Kac formula \eqref{FCTR}. By  deriving it with respect to $x$, we obtain
\begin{equation}
\frac{\partial V}{\partial x}(t,x)=  \E \left[ \int_0^{t}  h'(X^{t,x}_s) \frac{\partial X^{t,x}_s }{\partial x} ds \right] + \E\left[ g'(X^{t,x}_{t})\frac{\partial X^{t,x}_{t} }{\partial x}\right];
\end{equation}
 \item in Section \ref{Sec:4} we prove that, by extending a result of \cite{DeuZa} to time inhomogeneous diffusions and coefficients which are $C^1$ in the $x$ variable outside of a finite set, after the first hitting time of the boundary $\partial \shi = \{0,\zeta\}$ by the process $X_s$, then $\frac{\partial X^{t, x}_s }{\partial x} =0$;
 \item hence we can rewrite 
 \begin{equation}
  \frac{\partial V(t,x)}{\partial x}=  \E \left[ \int_0^{t}  h'(X^{t,x}_s) \frac{\partial X^{t,x}_s }{\partial x} 1_{\{ s< \tau^{t,x,X}_\partial \}}  ds\right]  + \E\left[ g'(X^{t,x}_{t})\frac{\partial X^{t,x}_{t} }{\partial x}1_{\{ t< \tau^{t,x,X}_\partial \}} \right];
 \end{equation}
where $\tau^{t,x,X}_\partial$ is the first hitting time of the boundary $\partial \shi $ by the process $X_s^{t,x}$;
\item since $ h'$ is bounded, in order to obtain a bound on $   \frac{\partial V(t,x)}{\partial x}$ it suffices to find a  bound on $\frac{\partial X^{t, x}_s }{\partial x}$ and  on $\P(s< \tau^{t,x,X}_\partial)$.
\end{enumerate}
We obtain a majorant to $\P(s< \tau^{t,x,X}_\partial)$ in the following way. 
In the case of arithmetic Brownian motion reflected at $0$ (resp. $\zeta$), a closed form spectral expansion is available for the law of the first hitting time to $\zeta$ (resp. $0$), see  \cite[Prop. 2]{Lin} but, since our process has non constant coefficients, we first perform a time change to obtain one with constant diffusion, $Z$, see \eqref{Z_def}. Such time change modifies the drift of the new process from $b(t-s,x)$ into $\hat b(t-s,x) = \frac{b(t-s,x)}{a(x)^2}$. At this point, the law of the first hitting time of $Z$ to  $\partial \shi = \{0,\zeta \}$ is bounded from above by the law of the first hitting time to $0$ (resp. $\zeta$) of a arithmetic Brownian motion reflected at $\zeta$ (resp. $0$), with drift $\mu_u =  \sup_{s \in [0,t],x\in \bar \shi} \hat b (s,x)$ (resp. $\mu_l = \inf_{s \in [0,t],x\in \bar \shi} \hat b (s,x)$), see \eqref{bound_tauXu} (resp. \eqref{bound_tauXd}).


By so doing, we obtain two  bounds that hold simultaneously on $\left| \frac{\partial V}{\partial x} \right|$. The first one, see Theorem \ref{316_bound_a}, is more  accurate around $x = 0$  while  the second one, proved in Theorem \ref{317_bound},  is accurate near  $x= \zeta$. 


By taking the minimum of these two bounds, we establish the main result, Theorem \ref{final_thm} below.  The constants $c_a^u, c^d_a$  below are introduced in  Hypothesis \ref{hyp}. 
Moreover, the functions $\shq_0$ (resp. $\shq_1$), defined in Theorems \ref{316_bound_a}  (resp. \ref{317_bound}),  only depend on time and on the coefficients of the problem, and are given in closed form. Each of them decomposes into two contributions, denoted by  $\shx$ and $\shy$ (see \eqref{shx} and \eqref{shy}). The first term, $\shx$, captures the contribution of the first eigenvalue of the spectral expansion of $\P(s< \tau^{\mu_u,x,S}_0)$ and $\P(s< \tau^{\mu_l,x,S}_\zeta)$,  while $\shy$ accounts for the tail contributions.


\begin{thm}\label{final_thm}
Assume Hypothesis \ref{hyp} is verified. Then, for any $t \geq 0$ and $x\in  \bar \shi$, the following holds.
\begin{equation}  \label{eq_380a_f}  
\left| \frac{\partial V}{\partial x}(t,x) \right|\leq 2 \frac{c^u_a }{ c^d_a } \min \left\{ \frac{x}{\zeta} \exp\left(\frac{-\mu_u}{\sigma^2} x\right) \shq_0(t,\mu_u,\sigma) , \frac{\zeta-x}{\zeta} \exp\left(-\frac{\mu_l}{\sigma^2} (x-\zeta)\right)   \shq_1(t,\mu_l,\sigma) \right\}
\end{equation}

\end{thm}
 
The dependence on  $x$  in  \eqref{eq_380a_f} originates from a uniform bound  on the coefficients $c_n(x)$, established in Proposition~\ref{prop_cn_bound}. These coefficients, introduced in \eqref{tauBound}, appear in the spectral expansion of $\P(s<\tau^{t,x,X}_0)$ and $\P(s<\tau^{t,x,X}_\zeta)$.

Additionally,  Corollary \ref{lim_infty} provides a sufficient condition ensuring that $\limsup_{t\ra \infty} \left| \frac{\partial V}{\partial x}(t,x) \right|$ remains bounded. Notably, this condition depends only on the sign of $\shd(\lambda_1,\sigma)$, see \eqref{D_def_a}, which is a linear function of  the first eigenvalue of the  spectral expansion.

Finally, we refine existing results on the derivative of the solution to a reflected SDE with respect to its initial condition,  $\frac{\partial X}{\partial x}$,   in the one dimensional case. In higher dimensions, \cite{DeuZa} proved that the map $x\mapsto \frac{\partial X}{\partial x}$ is continuous under the assumption of a  $C^1$ drift and constant diffusion. 
More recently, \cite{crisan2023} established a probabilistic representation for the derivative of the semigroup corresponding to a diffusion process killed at the boundary of a half interval, assuming time independent coefficients of class $C^2$. 

In the present work, we show that $x\ra \frac{\partial X}{\partial x}$ is continuous under weaker regularity requirements: specifically, the drift (resp. diffusion) is only required to be $C^1$ (resp. $C^2$) outside of a finite set, and the drift may also depend on time.

The main motivation of this work  stems from optimal control theory. When dealing with optimal control problems, a \textit{candidate} optimal control is often available but proving that it is actually optimal is often involved when the HJB does not admit a closed form solution. Indeed, it is necessary to prove that the candidate value function $V$ satisfies an inequality when the candidate optimal control is replaced by any other control, see for example \cite[Thm. 11.2.2]{Oksendal2003} and \cite[Thm. 8.1]{fleming2006}. In order to prove such inequality, explicit estimates on the derivatives of $V$ are often necessary when  a closed form solution is not available. 
This work will be applied to prove that the guessed optimal strategy of the sailboat trajectory optimization problem formulated in \cite{CDV1} is indeed optimal.

It is well known that the evaluation of option Greeks is a computationally expensive problem. Usually finite difference methods, Monte Carlo simulations (see \cite{CassFriz2006}), or perturbation analysis (see \cite{Glasserman1991}) is employed. The methodology we present here  is amenable to provide bounds for many Greeks with very little computational cost.
Indeed, the only quantity to be computed is the first eigenvalues of the Sturm-Liouville problem (see beginning of Section \ref{Sec:3}), which has a closed form expression, see \eqref{a_def_eq}.

Bounds on the gradient of the solution to a parabolic partial differential equation have been studied in various forms. The bounding function is usually a constant or a function of time alone. To the best of our knowledge, no results are available where the bound is expressed as a function of both space and time. These works typically consider either Dirichlet boundary conditions or PDEs defined on the whole $\R^n$. 

The first type of results concerns pointwise bounds on the gradient of the solution. 
 For instance, \cite{BartierSouplet2004} establishes, for semilinear PDE with Dirichlet boundary conditions, a global bound of the form $|\nabla_x V(t,x)| \leq C$, where $V$  is the solution to the parabolic PDE and  $C$ is a constant that is not provided in explicit form. More recently, \cite{Barles2021} extended the  Bernstein method to obtain local gradient estimates for solutions to second-order nonlinear elliptic and parabolic equations. For linear equations, similar results  appear in  \cite{LeBalch2021}, where the constant $C$ is expressed in the form of an exponential proportional to the $L_\infty$ norms of the coefficients, the diameter of the boundary, and the time. In a note, it is written that they also hold for homogeneous Neumann boundary conditions and null potential.
In the case of linear equations with constant coefficients in $\R^n$,  \cite{Kresin2020} derives a sharp bound of the form $|\nabla_x V(t,x)| \leq C(t) \|V(0,\cdot)\|_{L_p}$ and, in addition, also provides an explicit expression for the constant $C$. The bound is sharp, in the sense that $C$ cannot be reduced, and depends explicitly on time. 
Pointwise estimates for a class of parabolic PDEs with potential terms are also provided in \cite{Montoro2018}.

\cite{Han1998} and \cite{Kukavica2024} present pointwise Schauder estimates for the derivatives of a parabolic PDE in $\R^n$ that depend on the coefficient and non-homogeneous terms at this particular point

Another line of research concerns gradient bounds in  $L_p$ or  H\"older spaces. 
In \cite{Engler1990}, $L_p$ gradient estimates are proved for scalar  quasilinear and semilinear parabolic PDEs, and a special case of linear PDE in $\R^n$, where the scalar equations are coupled only through the non-homogeneous term. Bounds are of the form $\|\nabla_x V(t,\cdot)\|_{L_p} \leq \exp(C t)\|\nabla_x V(0,\cdot)\|_{L_p}$, exhibiting an explicit dependence on time. This work also addresses the long-time behavior of the gradient, showing that under certain conditions on the coefficients, the gradient remains bounded as $t \ra \infty$.

Results for Neumann boundary conditions are more limited. Among them, \cite{Lunardi} establishes Schauder estimates (see also Proposition \eqref{existSol} below). These results ensure Hölder continuity of the derivatives and, in combination with the boundary conditions, yield pointwise bounds in the interior of the domain. However, while the constant $C$ is known to depend on the coefficients, it is not given in explicit form. Classical references for Schauder estimates in this setting include \cite{LUS}, \cite{Fri}, \cite{solonnikov2016}. 

Another possible route is to differentiate the PDE, so that the new unknown $\frac{\partial V}{\partial x}$ satisfies a parabolic PDE with Dirichlet boundary conditions. By so doing, it is possible to apply results on pointwise bounds to the solution of a Dirichlet problem. However, this approach requires  at least $C^{1 +\alpha}$ regularity on the coefficients for some $\alpha >0$.


Despite this extensive literature, we are not aware of any work that provides  gradient bounds which are pointwise, with explicit closed-form expressions for the constants, precisely the contribution of the present work.

The paper is organized as follows. 
 
In Section \ref{Sec:2} we state  Hypothesis \ref{hyp} that will hold throughout the article and we derive a Feynman-Kac formula for $V$ in Proposition \ref{FeyC}. 
The main result of Section \ref{Sec:3} is Proposition \ref{Rem:3.4} which proves an upper bound on the law of the first hitting time to the boundary of the process \eqref{FC}.
In Section \ref{Sec:4} we prove that the derivative of the solution to \eqref{FC} is continuously differentiable, see Theorem \ref{ContBnd}.
The main Theorem \ref{final_thm} of the work is proved in Section \ref{Sec:5} together with some intermediary results. Finally, Corollary \ref{lim_infty} exhibits a sufficient condition under which $\limsup_{t\ra \infty}\left|\frac{\partial V}{\partial x}\right| <\infty$.

\section{Hypotheses and Feynman-Kac formula}\label{Sec:2}

We introduce  the main hypotheses on the coefficients and we prove a Feynman-Kac representation formula which will be used in the sequel to derive the estimates. 

For $\a \in ]0,1[$, we recall the parabolic H\"older spaces and norms $ C^{\alpha}$, $ C^{\alpha/2, \alpha}$, $C^{1+\alpha/2, 2+\alpha}$ defined, for example, in \cite{Lunardi}, Chapter 5, Section 1. Also, in the sequel, $\R_+ = \{t\in \R:t\geq 0 \}$.

In literature, usual regularity assumptions on the coefficients of the PDE are of H\"older or Lipschitz type. In the present work, in order to prove existence and continuity on the $x$-derivative of the solution to the SDE \eqref{FC}, we do need to require at least $C^1$ (resp. $C^2$) regularity of $h,b$ (resp. $a$)  outside of a finite set  because, in order to prove our result, we will need to differentiate those  coefficients.  
Throughout the paper, we assume that the coefficients of the PDE satisfy the following conditions.
\begin{hyp}\label{hyp}
The functions $a: \bar \shi \mapsto \R_+$, $b,h: \R_+ \times \bar \shi \mapsto \R$, $g: \bar \shi \mapsto \R$ are bounded and measurable.
Moreover, there exists $\a \in ]0,1[$ and $c, c_a^d,  c_a^u >0$ and a finite set $ \shg = \cup_{i=1}^g \hat x _i\subset \shi$ such that for any $t,t_1,t_2 \in \R_+$, $x_1,x_2\in \bar \shi$:
\begin{itemize}
\item[1)] $h \in C^{\alpha /2,\alpha}(\R_+\times \bar \shi)\cap  C^{\alpha /2,1}(\R_+\times \bar \shi \setminus \shg)$, $g \in C^{2+\alpha}(\bar \shi)$, 
\item[2)] $a\in C^1(\bar \shi) \cap C^2(\bar \shi \setminus \shg)$,
\item[3)]  $\left| b(t,x_1) - b(t,x_2) \right| \leq c  \left|x_1-x_2 \right|$ and $b \in C^{\alpha /2,1}(\R_+ \times  \bar \shi \setminus \shg)$,
\item[4)]  $c_a^d \leq a(x) \leq c_a^u$ for every $x \in  \bar \shi$,
\item [5)] $x\ra \frac{d^2 a}{dx^2}$ and $x \ra \frac{\partial b(t,x)}{\partial x}$ are bounded on $\bar \shi \setminus \shg$,
\item [6)] $\frac{d g(x)}{dx}|_{x=0} = 0$,  $\frac{d g(x)}{dx}|_{x=\zeta} = 0$.
\end{itemize}
\end{hyp}

\begin{rem}
With non-homogeneous boundary conditions, in order to apply Remark \ref{chg_var}, we need the following supplementary hypothesis: $n, m \in  C^{\alpha/2 +1}(\R_+)$ and Item 6) should be replaced with $\frac{d g(x)}{dx}|_{x=0} = n(0)$,  $\frac{d g(x)}{dx}|_{x=\zeta} = m(0)$.
 $n,m$ are more regular that what is strictly required for existence and uniqueness. This is because, when applying the change of variables of Remark \ref{chg_var}, then their time derivative enters in the driver $h$ of the PDE.
\end{rem}


\begin{prop}\label{existSol}
Assume Hypothesis \ref{hyp}. Then, the Neumann problem \eqref{freeB2aa}, \eqref{freeB2cc}, \eqref{freeB2ee} has a unique $C^{1+\a/2,2+\a}([0,T] \times \bar \shi)$ solution, for any $T>0$.
Moreover, there exists a constant $C$ such that
\begin{equation}\label{bound_sol}
\| V\|_{C^{1+\alpha/2, 2+\alpha}([0,T]\times \bar \shi)} \leq C \left( \|h \|_{C^{\alpha/2,\alpha }([0,T]\times \bar \shi)} +\|g \|_{C^{2+\alpha }(\bar{\shi})} \right).
\end{equation}
\end{prop}

\begin{proof}
Fix $T>0$.
The coefficient $b$ is Lipschitz in $x$ by Hypothesis \ref{hyp}, hence H\"older continuous for any $\a  \in ]0,1[$. In particular, $b \in C^{\a/2,\a}([0,T] \times \bar \shi)$, $a \in C^{\a}( \bar \shi)$ and $h\in  C^{\alpha /2, \alpha}([0,T]\times \bar \shi)$. Condition 5.1.67 of \cite{Lunardi} is also  satisfied by Item 6).
The result then follows from Corollary 5.1.22 of \cite{Lunardi}.
\end{proof}


The Proposition below provides a stochastic representation to the solution of the PDE.

\begin{prop}\label{FeyC}
Assume Hypothesis \ref{hyp} and let $V$ be a solution to \eqref{freeB2aa}--\eqref{freeB2ee}. Then, for $t\geq 0$,
\begin{equation}\label{FCTR}
V(t,x) = \E \left[ g(X_{t}^{t,x}) \right] + \E \left[ \int_0^{t} h(t-s,X_s^{t,x}) ds\right].
\end{equation}
\end{prop}

\begin{proof}
Proposition \ref{existSol} entails that $V \in C^{1,2}([0,T] \times \bar \shi)$, for any $T>t$. By It\^{o}'s formula,
\begin{align}
d V(t-s,X_s^{t,x}) &=  -\frac{\partial V}{\partial s}(t-s,X_s^{t,x}) ds +\frac{\partial V}{\partial x}( t-s, X_s^{t,x}) b(t-s, X_s^{t,x}) ds\\
&\qquad  + \frac{\partial  V}{\partial x}(t-s,X_s^{t,x}) \s  a(X_s^{t,x}) dB_t \\
&\qquad  + \frac{\partial  V}{\partial x}(t-s,X_s^{t,x}) ( dL_s^{t,x} - dU_s^{t,x}) \label{1.42}\\
&\qquad + \frac{\sigma^2}{2} \frac{\partial^2  V}{\partial x^2}(t-s,X_s^{t,x})  a^2(X_s^{t,x}) ds. 
\end{align}
Since $ dL_s^{t,x}$, (resp. $ dU_s^{t,x}$) is positive only when the process hits $0$ (resp. $\zeta$) and here by the Neumann boundary condition $ \frac{\partial V}{\partial x}(s,0)= 0$ (resp. $ \frac{\partial V}{\partial x}(s,\zeta)=0$, the term \eqref{1.42} is always equal to zero.
Using \eqref{freeB2aa}, we can rewrite the above as 
\begin{align}
d V(t-s,X_s^{t,x}) = -h(t-s,X_s^{t,x})  + \frac{\partial  V}{\partial x}(t-s,X_s^{t,x}) \s a(X_s^{t,x})  dB_s.
\end{align}
Integrate from $0$ to $t$  and take expectation. The stochastic integral is a martingale since  $\frac{\partial V}{\partial x}(s,x)$ is uniformly continuous  on $[0,T]\times \bar \shi$ by Proposition \ref{existSol}, so we obtain
\begin{align}\label{eq26a}
 V(t,x) -\E \left[ V(0,X_{t}^{t,x}) \right] = \E \left[\int_{0}^{t}  h(t-s,X_s^{t,x}) ds \right].
\end{align}
Using the initial condition \eqref{freeB2ee}, the Proposition is proved.
\end{proof}


\section{A bound on hitting times probabilities.}\label{Sec:3}

In this section, we establish an upper bound on the distribution of the first hitting time to the boundary $\partial \shi $ of the process \eqref{FC}. Corollary \ref{Cor:1.43} demonstrates that this law is bounded by the law of  first hitting time to the boundary of a process $Z$ with  constant diffusion, which is derived from \eqref{FC} by means of a time change. In turn, Proposition \ref{Rem:3.4}, which is the main result of the Section, provides an upper bound for the law of the first hitting time of $Z$  in terms of that of a reflected arithmetic Brownian motion, and the latter can be expressed using a spectral expansion.

We begin by recalling a result about the  spectral expansion representation for the first hitting time of reflected arithmetic Brownian motion with drift. 

Let $\mu \in \R$, $s\geq 0$, $x\in \bar \shi$.  Define  $\tau^{\mu,x,S}_y:= \inf \{s\geq 0:\, S^{\mu,x}_s  = y \} $, for $y<x$ (resp. $y>x$), where $S^{\mu,x}$ is the $\P_x$-unique strong solution to the arithmetic Brownian motion reflected at $\zeta$ (resp. $0$): 
\begin{align}
dS^{\mu,x}_s &= \mu ds + \s dB_s - dU^{\mu,x}_s \label{eq:5.1a} \\
S^{\mu,x}_s & \leq \zeta,
\end{align}
respectively,
\begin{align}
dS^{\mu,x}_s &= \mu ds + \s dB_s + dL^{\mu,x}_s \label{eq:5.1b} \\
S^{\mu,x}_s & \geq 0.
\end{align}

By \cite{Lin}, Proposition 2 and Remark 3, it follows that 
\begin{equation}\label{tauBound}
\P (s\leq \tau^{\mu,x,S}_y < \infty) = -\sum_{n=1}^\infty \exp(-\lambda_{n,y} s) \frac{\psi(x,\lambda_{n,y}) }{\lambda_{n,y} \psi_\lambda (y,\lambda_{n,y}) } = \sum_{n=1}^\infty \exp(-\lambda_{n,y} s) c_n(x),
\end{equation}
where $\psi(x,\lambda)$ is the unique solution (up to a constant multiplication factor independent from $x$) to
\begin{align}\label{SL1}
&\um \s^2 \psi'' + \mu \psi' + \lambda \psi = 0, \qquad x >0,
\end{align}
with boundary condition $\psi'(\zeta) = 0$, (resp. $\psi'(0) = 0$).
The eigenvalues $\{\lambda_{n,y}\}_{n=1}^\infty$, $0<\lambda_{1,y}<\lambda_{2,y}<\dots \lambda_{n,y} \ra \infty$ as $n\ra \infty$ are the simple positive zeros of $\psi(0,\lambda)$ (resp. $\psi(\zeta,\lambda)$), that is the values such that $\psi(0,\lambda_{n,0}) = 0$  (resp. $\psi(\zeta,\lambda_{n,\zeta}) = 0$).

The solution to \eqref{SL1} and the $c_n$ are provided explicitly, in the case of Brownian motion with drift, in (3.16) of \cite{BMhit}. We recall below the main results.

\textit{Case 1: $\mu \geq 0$, hitting time down ($y<x$).}

Here, $S^{\mu,x}$ is reflected at $\zeta$ and the eigenvalues are the simple positive zeros of \eqref{SL1} with boundary condition $\psi'(\zeta) = 0$. We denote the eigenvalues by $\lambda^{u,\mu}_{n,y}$, where the $u$ in the superscript is meant to distinguish it from the $\lambda^{l,\mu}_{n,y}$ defined below related to the \textit{hitting time up problem} and $\mu$ is meant to keep track of the drift of the process.

For $n=1$, $\mu >0$, the expression of  $\lambda_{1,y}$ depends on the value of the ratio $\frac{\s^2}{\mu(\zeta-y)}$, namely
\begin{equation}\label{lambda_def}
\lambda^{u,\mu}_{1,y} = 
\left\{ \arraycolsep=1.2pt\def\arraystretch{1.8}
\begin{array}{rll}
&\frac{\mu^2 - v_{\mu,y}^2}{2\s^2}, &\frac{\s^2}{\mu(\zeta-y)} \in ]0,1[,\\
&\frac{\mu^2}{2\s^2},  &\frac{\s^2}{\mu(\zeta-y)} = 1, \qquad \text{and} \, \lambda^{u,\mu}_{n,y} = \frac{(\alpha^{u,\mu}_{n,y})^2\s^2}{2(\zeta-y)^2} + \frac{\mu^2}{2\s^2}\\
& \frac{(\alpha^{u,\mu}_{1,y})^2\s^2}{2(\zeta-y)^2} + \frac{\mu^2}{2\s^2}, &\frac{\s^2}{\mu(\zeta-y)} \in ]1,+\infty[,
\end{array}
\right.
\end{equation}
where $\a^{u,\mu}_{n,y}$ is the unique solution, in the interval $[(n-1)\pi, (n-1)\pi + \pim[$, to 
\begin{equation}\label{a_def_eq}
\tan \alpha = \frac{\s^2}{\mu(\zeta-y)}  \alpha.
\end{equation}
The $v_{\mu,y}$ appearing in \eqref{lambda_def} is the unique solution in $]0,\mu[$ to 
\begin{equation}\label{x0_def}
 \frac{\s^2}{2(\zeta-y)} \ln \left(\frac{\mu+v}{\mu-v} \right) = v,
\end{equation}
which exists if $\frac{\s^2}{\mu(\zeta-y)} \in ]0,1[$.
The $c_{n}(x)$ introduced in \eqref{tauBound} have the following form
\begin{equation}\label{cnDef1}
c_{n,y}^{u,\mu}(x) = \frac{\sqrt{\mu^2-2\lambda^{u,\mu}_{n,y} \s^2}}{\mu - 2 \lambda^{u,\mu}_{n,y} (\zeta-y)}\left[
\exp\left(\frac{-\mu^-(\lambda^{u,\mu}_{n,y})}{\s^2}(x-y)\right) - \exp \left(\frac{-\mu^+(\lambda^{u,\mu}_{n,y})}{\s^2}(x-y)\right)\right],
\end{equation}
where 
\begin{equation}\label{mu_pm_def}
\mu ^{\pm}(\lambda) = \mu \pm \sqrt{\mu^2-2\lambda \s^2}. 
\end{equation}

The  case where $\mu = 0$ can be tackled explicitly observing that $\psi(x,\lambda) = \cos \left((\zeta-x) \sqrt{\frac{2 \lambda^{u,0}_{n,y}}{\s^2}}\right)$:
\begin{align}\label{mu0_case}
&\lambda^{u,0}_{n,y} = \frac{(2n-1)^2 \s^2\pi^2}{8(\zeta-y)^2},\\
& c_{n,y}^{u,0}(x) = \frac{(-1)^{n+1}4}{(2n-1)\pi} \sin \left( \frac{(2n-1)\pi}{2}\frac{(x-y)}{(\zeta-y)}\right). 
\end{align}
It can also  be seen as the limiting case $\frac{\s^2}{\mu(\zeta-y)} \ra \infty$.

\textit{Case 2: $\mu \geq 0$, hitting time up ($y>x$).}

If $\mu >0$, the process $S^{\mu,x}$ is reflected at $0$ instead of $\zeta$. 
We remark that in this case there are no eigenvalues in $[0, \frac{\mu^2}{2\s^2}[$, so  there aren't the three cases as in \eqref{lambda_def} and all the eigenvalues have the same expression \eqref{lambda_def1}.
Equation (3.23) of \cite{BMhit} reads as
\begin{equation}\label{lambda_def1}
\lambda^{l,\mu}_{n,y} = \frac{(\alpha^{l,\mu}_{n,y})^2\s^2}{2y^2} + \frac{\mu^2}{2\s^2},
\end{equation}
where $\a^{l,\mu}_{n,y}$ is the unique solution, in the interval $]\pim +(n-1)\pi, \pi +(n-1)\pi]$ for $n\geq 1$, to 
\begin{equation}\label{a_def_eq_1}
\tan \alpha^{l,\mu} = -\frac{\s^2}{\mu y}  \alpha^{l,\mu},
\end{equation}
and
\begin{equation}\label{cnDef2}
c_{n,y}^{l,\mu}(x) = \frac{\sqrt{\mu^2-2\lambda^{l,\mu}_{n,y} \s^2}}{\mu + 2 \lambda^{l,\mu}_{n,y} y}\left[
\exp\left(\frac{-\mu^-(\lambda^{l,\mu}_{n,y})}{\s^2}(x-y)\right) - \exp\left(\frac{-\mu^+(\lambda^{l,\mu}_{n,y})}{\s^2}(x-y)\right)\right].
\end{equation}
The case $\mu = 0$ can be handled explicitly as well observing that $\psi(x,\lambda) = \cos \left(x\sqrt{\frac{2 \lambda}{\s^2}} \right)$:
\begin{align}\label{mu0_case_1}
&\lambda^{l,0}_{n,y} = \frac{(2n-1)^2 \s^2\pi^2}{8y^2},\\
& c_{n,y}^{l,0}(x) = \frac{(-1)^{n+1}4}{(2n-1)\pi} \sin \left( \frac{(2n-1)\pi}{2}\frac{(y-x)}{y}\right). 
\end{align}

\begin{rem}\label{mu_neg_rem}
If $\mu < 0$ the result follows by a reflection.
\begin{enumerate}
\item \textit{Hitting down}. If $y \in [0,x[$, then we can do a reflection along the axis $\zeta/2$, hence $y$ and $x$ are replaced by $\zeta -y$ and $\zeta-x$ and the drift becomes $-\mu >0$.  Now apply \textit{Case 2} with these parameters.
\item \textit{Hitting up}. If $y \in ]x,\zeta]$, then we can do again a reflection along the axis $\zeta/2$ and apply \textit{Case 1} with the new parameters.
\end{enumerate}
\end{rem}

Now, a useful Lemma that will be applied in the next Proposition \ref{prop_cn_bound} to bound  $c^{u,\mu}_{n,y}$ defined in \eqref{cnDef1}.
\begin{lem}\label{sign_cn}
Assume Hypothesis \ref{hyp} is in force. Let $(\lambda^{u,\mu}_{n,y})_{n\geq 1}$ be the sequence of eigenvalues defined in \eqref{lambda_def}. For any $\mu > 0$, $\sigma>0$, $y\in [0,\zeta[$, then
\begin{equation}\label{eq_280_sign}
\mu - 2 \lambda^{u,\mu}_{n,y}(\zeta-y) \leq 0, \qquad n\geq1.
\end{equation}
In particular, for $\frac{\s^2}{\mu(\zeta-y)}  = 1$ and $n = 1$, \eqref{eq_280_sign} is equal to zero. Nonetheless, in this case, $c_{1,y}^{u,\mu}$ in \eqref{cnDef1} is still well-defined.
\end{lem}
\begin{proof}
Below, we denote $\lambda_{n} = \lambda^{u,\mu}_{n,y}$.
Since $\lambda_{n}$ is an increasing sequence, it is enough to prove the claim for $n=1$.

\textit{Case 1.} $\frac{\s^2}{\mu(\zeta-y)} <1$. Then $\lambda_{1} = \frac{\mu^2 - v_{\mu,y}^2}{2\s^2}$, hence
\begin{align}\label{eq_230}
\mu - 2 \lambda_{1}(\zeta-y) = \mu - \frac{\mu^2 -  v_{\mu,y}^2}{\s^2}(\zeta-y) \geq \mu - \frac{\mu^2 (\zeta-y)}{\s^2}.
\end{align}
Dividing \eqref{eq_230} by $\mu$, we obtain
$1 - \frac{\mu(\zeta-y) }{\s^2}$, which is negative.

\textit{Case 2.} $\frac{\s^2}{\mu(\zeta-y)} = 1$.
Recall that $\lambda_{1} = \frac{\mu^2}{2\s^2}$ according to \eqref{lambda_def}, hence
\begin{align}
\mu - 2 \lambda_{1}(\zeta-y) = \mu -  \frac{\mu^2}{\s^2}(\zeta-y) =0.
\end{align}
We now prove that $c_{1,y}^{u,\mu}$ is well-defined. We compute $$-\lim_{\lambda \searrow \lambda_{1}}   \frac{\psi(x,\lambda) }{\lambda \psi_\lambda (y,\lambda) }  =\lim_{\lambda \searrow \lambda_{1}}  c_{1,y}^{u,\mu}(\lambda).$$ 
Let $\lambda$ converge from above to $\lambda_1$. By \eqref{cnDef1},
\begin{align}\label{cnDef10}
\lim_{\lambda \searrow \lambda_{1}} c_1(\lambda) &= \lim_{\lambda \searrow \lambda_{1}} -2\frac{\sqrt{2\lambda \s^2-\mu^2}}{\mu - 2 \lambda (\zeta-y)} \exp\left(\frac{-\mu}{\s^2}(x-y)\right)  \sin\left(\frac{\sqrt{2\lambda \s^2- \mu^2}}{\s^2} (x-y)\right)\\
& = \lim_{\lambda \searrow \lambda_{1}} -2\frac{2\lambda \s^2-\mu^2}{\mu - 2 \lambda (\zeta-y)} \exp\left(\frac{-\mu}{\s^2}(x-y)\right) \frac{x-y}{\s^2} ,
\end{align}
where in the last row, since $\sqrt{2\lambda \s^2- \mu^2}  \ra 0 $ as $\lambda \ra \lambda_1$, we used $ \sin\left(\frac{\sqrt{2\lambda \s^2- \mu^2}}{\s^2} (x-y)\right) \approx  \frac{\sqrt{2\lambda \s^2- \mu^2}}{\s^2} (x-y)$.
Using the relation $\frac{\s^2}{\mu(\zeta-y)} = 1$, the above becomes
\begin{align}\label{cnDef11}
 \lim_{\lambda \searrow \lambda_{1}} -2\frac{2\lambda \s^2-\mu^2}{\mu - 2 \lambda \frac{\s^2}{\mu}} \exp\left(\frac{-\mu}{\s^2}(x-y)\right) \frac{x-y}{\s^2}= 2\mu \exp\left(\frac{-\mu}{\s^2}(x-y)\right)\frac{x-y}{\s^2} .
\end{align}
After repeating the above argument for $\lambda \nearrow \lambda_{1}$, the claim is proved.

\textit{Case 3.} $\frac{\s^2}{\mu(\zeta-y)} > 1$.
For $n \geq 1$
\begin{align}\label{lambda_n_3}
\mu - 2 \lambda_{1}(\zeta-y) = \mu - \frac{(\alpha^{u,\mu}_{1,y})^2\s^2}{(\zeta-y)} - \frac{\mu^2}{\s^2}(\zeta-y).
\end{align}
Recalling that $\alpha_1 := \alpha^{u,\mu}_{1}$ is the unique solution in $[0,\pim[$ of the equation \eqref{a_def_eq}, we see that 
\begin{equation}\label{a_def_eq11}
\mu = \frac{\alpha_1}{\tan(\alpha_1)}\frac{\s^2}{\zeta-y}.
\end{equation}
Replacing \eqref{a_def_eq11} in \eqref{lambda_n_3}, we obtain
\begin{align}
\frac{\s^2}{\zeta-y}\left(\frac{\alpha_1}{\tan(\alpha_1)} -\alpha^2_1 - \frac{\alpha_1^2}{\tan^2(\alpha_1)} \right).
\end{align}
Since $\alpha_1 \in [0,\pim[$, if $f(\alpha):= \left(\frac{1}{\tan(\alpha)} -\alpha - \frac{\alpha}{\tan^2(\alpha)}\right) \leq 0$ for every $\a \in [0,\pim[$, the claim would be proved.
We first observe that $\lim_{\alpha \ra 0}f(\alpha) = 0$. We now prove that $\frac{d f}{d\alpha} \leq 0$. After some straightforward calculations, the following holds  
\begin{equation}
\frac{d f}{d\alpha} = 2 \left(  \frac{\alpha}{\tan(\alpha)} - 1\right),
\end{equation}
which is non-positive for $\a \in [0,\pim[$ and the Lemma is proved.
\end{proof}

The proposition below proves  bounds on  $ c_{n,y}^{u,\mu}(x) $ and $ c_{n,y}^{l,\mu}(x) $ that will be used in the proofs of Theorem \ref{316_bound_a} and \ref{317_bound}.

\begin{prop}\label{prop_cn_bound}
Let $\mu \geq 0$, then
for any $x\in \bar \shi$  the following holds.
\begin{enumerate}[(1)]
 \item If  $\frac{\sigma^2}{\mu (\zeta-y)} \geq 1$,  and $ y\in [0,x[$, then 
  for $n\geq 1$
   \begin{equation}
\left| c_{n,y}^{u,\mu}(x)  \right| \leq   2  \frac{ x-y}{\zeta-y}  c_{\alpha,y}^{u,\mu} \exp\left(-\frac{\mu}{\sigma^2} (x-y)\right),
\end{equation}
  where 
  \begin{equation}
   c_{\alpha,y}^{u,\mu} = \frac{(\alpha^{u,\mu}_{1,y})^2 \sigma^2}{- \mu (\zeta -y)  + (\alpha^{u,\mu}_{1,y})^2 \sigma^2 + \frac{\mu^2 (\zeta-y)^2}{\sigma^2}};
   \end{equation}

 \item   if  $\frac{\sigma^2}{\mu (\zeta-y)} <1$ and $ y\in [0,x[$ then,
 \begin{enumerate}[(a)]
\item  for $n>1$,
  \begin{equation}
\left| c_{n,y}^{u,\mu}(x)  \right| \leq   2  \frac{ x-y}{\zeta-y} \exp\left(-\frac{\mu}{\sigma^2} (x-y)\right);
\end{equation}
\item for $n=1$,
\begin{equation}
\left| c^{u,\mu}_{1,y}(x) \right|  \leq 2 \frac{x-y}{\zeta-y} \exp\left(-\frac{\mu}{\sigma^2} (x-y)\right) \sinh\left(\frac{v_{\mu,y}}{\s^2}(\zeta-y)\right) \frac{\s^2 v_{\mu,y}}{\s^2 \mu -(\mu^2 -v_{\mu,y}^2)( \zeta-y)},
\end{equation}
where $v_{\mu,y}$ is defined in \eqref{x0_def};
\end{enumerate}
\item if $ y\in ]x,\zeta]$, then  for $n\geq1$
   \begin{equation}
\left| c_{n,y}^{l,\mu}(x)  \right| \leq   2  \frac{ y-x}{y} \exp\left(-\frac{\mu}{\sigma^2} (x-y)\right).
\end{equation}
 \end{enumerate}  
\end{prop}

\begin{proof}
\begin{enumerate}[(1)]
\item
Define $\a_n:= \alpha^{u,\mu}_{n,y}$ and $\lambda_n :=\lambda^{u,\mu}_{n,y}$, where $\alpha^{u,\mu}_{n,y}$ and $ \lambda^{u,\mu}_{n,y}$ are introduced, respectively, in \eqref{a_def_eq} and in \eqref{lambda_def}.
By \eqref{lambda_def} and \eqref{mu_pm_def} it follows that
\begin{equation}
\mu^{\pm}(\lambda_n) = \mu \pm i  \frac{\alpha_n \sigma^2}{\zeta -y},
\end{equation}
hence
\begin{equation}\label{eq3444}
\exp\left(\frac{-\mu^-(\lambda_n)}{\s^2}(x-y)\right) - \exp\left(\frac{-\mu^+(\lambda_n)}{\s^2}(x-y)\right) = \exp\left(-\frac{\mu}{\sigma^2} (x-y)\right) 2i \sin\left(\alpha_n\frac{x-y }{\zeta-y}\right).
\end{equation}
Additionally, 
\begin{equation}\label{eq3454}
\frac{\sqrt{\mu^2-2\lambda_n \s^2}}{\mu - 2 \lambda_n (\zeta-y)} = i \frac{\alpha_n \sigma^2}{\mu (\zeta-y)  - \alpha_n^2 \sigma^2 - \frac{\mu^2 (\zeta-y)^2}{\sigma^2}}.
\end{equation}
Replacing \eqref{eq3444} and \eqref{eq3454} in \eqref{cnDef1}, 
\begin{equation}\label{eq_3904}
 c_{n,y}^{u,\mu}(x)   =2 \frac{\alpha_n \sigma^2}{-\mu (\zeta-y)  + \alpha_n^2 \sigma^2 + \frac{\mu^2 (\zeta-y)^2}{\sigma^2}} \exp\left(-\frac{\mu}{\sigma^2} (x-y)\right) \sin \left(\alpha_n\frac{x-y }{\zeta-y}\right).
\end{equation}
We recall that, if $\mu = 0$, then \eqref{eq_3904} becomes \eqref{mu0_case}. Additionally,
observe that for the limit case $\frac{\s^2}{\mu_u (\zeta-y)} = 1$, \eqref{eq_3904} is still well-defined by \textit{Case 2)} of Lemma \ref{sign_cn}.

\textit{Case 3)} of Lemma \ref{sign_cn} shows  that  $\mu -2\lambda_n (\zeta-y)  \leq 0$. Hence,  by  using  \eqref{lambda_def} and changing sign, we obtain $- \mu (\zeta -y)  + \alpha_n^2 \sigma^2 + \frac{\mu^2 (\zeta-y)^2}{\sigma^2} \geq  0$. Moreover, recalling the elementary inequality $|\sin(x)| \leq x$ for $x\geq0$ , we deduce
\begin{align}
|c_{n,y}^{u,\mu}(x)| 
&\leq 2  \frac{x-y}{\zeta-y} \exp\left(-\frac{\mu}{\sigma^2} (x-y)\right) \frac{\alpha_n^2 \sigma^2}{- \mu (\zeta-y)  + \alpha_n^2 \sigma^2 + \frac{\mu^2 (\zeta-y)^2}{\sigma^2}}  .\label{eq348a_24}
\end{align}
Since $\frac{\sigma^2}{\mu (\zeta-y)} \geq 1$, then $- \mu (\zeta-y)  +  \frac{\mu^2 (\zeta-y)^2}{\sigma^2} \leq 0$, hence for $n\geq 1$,  the function
\begin{equation}\label{eq348a_241}
\a \ra f(\alpha) :=\frac{\alpha^2 \sigma^2}{- \mu (\zeta-y)  + \alpha^2 \sigma^2 + \frac{\mu^2 (\zeta-y)^2}{\sigma^2}}
\end{equation}
is monotonically decreasing.
Therefore \eqref{eq348a_24} can be bounded from above by
\begin{equation}\label{eq_348_34}
2  \frac{x-y}{\zeta-y}   \exp\left(-\frac{\mu}{\sigma^2} (x-y)\right) \frac{\alpha_1^2 \sigma^2}{- \mu (\zeta-y)  + \alpha_1^2 \sigma^2 + \frac{\mu^2 (\zeta-y)^2}{\sigma^2}}.
\end{equation}

\item 
Define, as in Item (1), $\a_n:= \alpha^{u,\mu}_{n,y}$ and $\lambda_n :=\lambda^{u,\mu}_{n,y}$.
By \eqref{lambda_def}, the first eigenvalue is structurally different, hence we will consider first $n>1$.
Observe that \eqref{eq_3904} also holds in this case. 
Since $\frac{\sigma^2}{\mu (\zeta-y)} <1$, then $- \mu (\zeta-y) + \frac{\mu^2 (\zeta-y)^2}{\sigma^2} > 0$, hence 
\begin{align}\label{eq348a4}
 \left|c_{n,y}^{u,\mu}(x) \right| = &2 \exp\left(-\frac{\mu}{\sigma^2} (x-y)\right)  \left| \frac{\alpha_n \sigma^2}{- \mu (\zeta -y)  + \alpha_n^2 \sigma^2 + \frac{\mu^2 (\zeta-y)^2}{\sigma^2}}   \sin\left(\alpha_n \frac{ x-y}{\zeta-y} \right) \right| \\
&\leq 2\exp\left(-\frac{\mu}{\sigma^2} (x-y)\right) \frac{ x-y}{\zeta -y} .
\end{align}

For $n=1$, $\lambda_1 = \frac{\mu^2 -v_{\mu,y}^2}{2\sigma^2}$, and $c_{1,0}^{u,\mu}$ is as in \eqref{cnDef1}, hence:
\begin{align}\label{eq_360}
\frac{\sqrt{\mu^2-2\lambda_1 \sigma^2}}{\mu-2\lambda_1(\zeta-y)} = \frac{v_{\mu,y}}{\mu -\frac{\mu^2 -v_{\mu,y}^2}{\sigma^2} (\zeta-y)},
\end{align}
and, since $\mu^\pm(\lambda_1 ) = \mu \pm v_{\mu,y}$,
\begin{align}\label{eq_361}
&\exp\left(\frac{-\mu^-(\lambda_1)}{\s^2}(x-y)\right) - \exp\left(\frac{-\mu^+(\lambda_1)}{\s^2}(x-y)\right)\\
&\qquad  = \exp\left(-\frac{\mu}{\s^2}(x-y)\right) \left( \exp\left(\frac{v_{\mu,y}}{\s^2}(x-y)\right) -\exp\left(-\frac{v_{\mu,y}}{\s^2}(x-y)\right)\right) \\
&\qquad=  2 \exp\left(-\frac{\mu}{\s^2}(x-y) \right) \sinh\left(\frac{v_{\mu,y}}{\s^2}(x-y)\right)\\
 &\qquad \leq 2 \exp\left(-\frac{\mu}{\s^2}(x-y)\right) \frac{ \sinh\left(\frac{v_{\mu,y}}{\s^2}(\zeta-y)\right)}{\zeta-y} (x-y).
\end{align}
Therefore, for $n=1$, by \eqref{eq_360} and \eqref{eq_361}, the conclusion follows.

\item 
Define $\a_n:= \alpha^{l,\mu}_{n,y}$ and $\lambda_n :=\lambda^{l,\mu}_{n,y}$, where $\alpha^{l,\mu}_{n,y}$ and $ \lambda^{l,\mu}_{n,y}= \frac{\alpha_n^2\s^2}{2 y^2} + \frac{\mu^2}{2\s^2}$ are introduced, respectively, in \eqref{a_def_eq_1} and in \eqref{lambda_def1}.

By \eqref{mu_pm_def},
\begin{equation}\label{eq_324a}
\mu^{\pm}(\lambda_n) = \mu \pm i  \frac{\alpha_n \sigma^2}{y},
\end{equation}
Using \eqref{eq_324a} in \eqref{cnDef2}, and arguing as in Item (1), we obtain
\begin{align}
&\left|c_{n,y}^{l,\mu}(x)\right| =\left| \frac{\sqrt{\mu^2-2\lambda_n \sigma^2}}{\mu+2\lambda_n y} \left(\exp\left(\frac{-\mu^-(\lambda_n)}{\s^2}(x-y)\right) - \exp\left(\frac{-\mu^+(\lambda_n)}{\s^2}(x-y)\right)\right)\right| \\ \label{eq_390_14}
&\qquad  = \left|2 \exp\left(-\frac{\mu}{\sigma^2} (x-y)\right)  \frac{\alpha_n \sigma^2}{-\mu y  - (\alpha_n)^2 \sigma^2 - \frac{\mu^2 y^2}{\sigma^2}}  \sin\left(\alpha_n\frac{ x-y}{y} \right)\right| \\
&\qquad \leq 2  \exp\left(-\frac{\mu}{\sigma^2} (x-y)\right)    \frac{ y-x}{y}.\label{eq348a_2_14}
\end{align}
Also here we recall that, if $\mu =0$, then \eqref{eq_390_14} becomes \eqref{mu0_case_1}.
\end{enumerate}
\end{proof}

\begin{rem}\label{rem_34}
For $n$ large enough, if the conditions of Item (1) of Proposition \ref{prop_cn_bound} are verified, the following bound holds.

Let $x \in \bar{\shi}$ and $y\in [0,x[$.  If $\mu >0$ and $\frac{\sigma^2}{\mu (\zeta-y)} \geq 1$ then, for any $\ve>0$, there exists $\bar n$ such that, for any $n \geq \bar n$,
  \begin{equation}
\left| c_{n,y}^{u,\mu}(x)  \right| \leq  2 (1+\ve)   \frac{ x-y}{\zeta-y} \exp\left(-\frac{\mu}{\sigma^2} (x-y)\right).
\end{equation}
Indeed, the function $f(\alpha)$, defined in  \eqref{eq348a_241}, is decreasing in $\alpha$ and $\lim_{\alpha \ra \infty} f(\alpha) = 1$.
The conclusion then follows recalling that, by \eqref{a_def_eq}, $\alpha^{u,\mu}_{n,y} \geq (n-1)\pi$.

\end{rem}

In the sequel, we will apply \eqref{tauBound} in the following cases:
\begin{enumerate}
\item drift $\mu \geq 0$, hitting time of the level $y=0$, starting from $x>0$, with reflection at $\zeta$: see \eqref{eq:5.1a}. Therefore we simplify notations and the solution to \eqref{a_def_eq}, \eqref{x0_def} and the eigenvalues defined in \eqref{lambda_def} read
\begin{align}
\alpha^{u,\mu}_{n} &:=\alpha^{u,\mu}_{n,y}; \label{alpha_u_def}\\
v_{\mu}&:= v_{\mu,0}; \label{v_mu_def}\\
\lambda^{u,\mu}_n &:= \lambda^{u,\mu}_{n,0} =  \frac{(\alpha^{u,\mu}_{n})^2\s^2}{2\zeta^2} + \frac{\mu^2}{2\s^2};\label{lambda_u_def}
\end{align}
\item drift $\mu > 0$, hitting time of the level $y=\zeta$ starting from $x<\zeta$, with reflection at $0$: see \eqref{eq:5.1b}. Therefore the solution to \eqref{a_def_eq_1} and the eigenvalues defined in \eqref{lambda_def1} read
\begin{align}
\alpha^{l,\mu}_{n} &:=\alpha^{l,\mu}_{n,y}; \label{alpha_l_def}\\
\lambda^{l,\mu}_n &:=\lambda^{l,\mu}_{n,\zeta} =  \frac{(\alpha^{l,\mu}_n)^2 \s^2}{2\zeta^2} + \frac{\mu^2}{2\s^2}.\label{lambda_l_def}
\end{align}
\end{enumerate}
The case $\mu < 0$ follows by a reflection as pointed out in Remark \ref{mu_neg_rem}. Also,  see \eqref{mu0_case} if $\mu= 0$.

The change-of-time formula below is useful to transform the process $(X^{t,x}_s)_{s\in [0,t]}$, defined in \eqref{FC}, into a process with a constant diffusion term $Z$, that is, to remove the function $a$. This is necessary in order to bound $Z$ from above and below with a arithmetic Brownian motion and apply the spectral expansion \eqref{tauBound}.

\begin{lem}\label{Lem:1.41}
Assume Hypothesis \ref{hyp}. Let $t\geq 0$ and $x\in \bar \shi$. Define
\begin{align}
A^{t,x}_s:=\int_0^s  a(X^{t,x}_r)^2dr,\qquad s \in [0,t],
\end{align}
with inverse
\begin{align}
A^{-1,t,x}_s:= \inf \{r \geq 0:\, A^{t,x}_r \geq s \}, \qquad s\in [0, A^{t,x}_{t}],
\end{align}
and $Z^{t,x}_s:=X^{t,x}_{A^{-1,t,x}_s}$. Then, for $s\in [0, A^{t,x}_{t}]$,
\begin{align}
&dZ^{t,x}_s= \sigma d\tilde{B}_s + \int_0^s b(t - A^{-1,t,x}_u, Z^{t,x}_u) \frac{1}{a(Z^{t,x}_u)^2}du +dL^{t,x,Z}_{s}-dU^{t,x,Z}_{s}, \label{Z_def}\\
&  Z^{t,x}_s\in \bar \shi \, a.s., \qquad Z^{t,x}_0=x,\\
	&1_{\{Z_s^{t,x} \neq 0 \}} dL^{t,x,Z}_{s}= 0,\qquad 1_{\{Z_t^{t,x} \neq \zeta \}} dU^{t,x,Z}_{s}= 0,
\end{align}
where $\tilde{B}_s$ is a $(\shf_s)$-Brownian motion.
\end{lem}

\begin{proof}
We will sometimes write $A^{-1}_s$ in place of $A^{-1,t,x}_s$ when this does not generate confusion.
By \eqref{FC}, for $s\in [0,t]$,
\begin{equation}
X^{t,x}_s = x + \sigma \int_0^s  a(X^{t,x}_r) dB_r +\int_0^s b(t - r, X^{t,x}_r) dr + L^{t,x}_s-U^{t,x}_s,
\end{equation}
hence, defining $Z^{t,x}_s:= X^{t,x}_{A^{-1}_s}$, it follows that
\begin{equation}\label{eq_36a}
Z^{t,x}_s=  x + \sigma  \int_0^{A^{-1}_s} a(X^{t,x}_r) dB_r +\int_0^{A^{-1}_s} b(t - r, X^{t,x}_r) dr  +L^{t,x}_{A^{-1}_s}-U^{t,x}_{A^{-1}_s},
\end{equation}
in particular $Z^{t,x}_s \in \bar \shi$ a.s. Let
\begin{align}
 L^{t,x,Z}_{s}:= L^{t,x}_{A^{-1}_s},\qquad U^{t,x,Z}_{s}:= U^{t,x}_{A^{-1}_s},
\end{align}
 then $	1_{\{Z_s^{t,x} \neq 0 \}} dL^{t,x,Z}_{s}= 0$, and $1_{\{Z_s^{t,x} \neq \zeta \}} dU^{t,x,Z}_{s}= 0$.
The stochastic integral can be rewritten as $\sigma \tilde{B}_s$ where $\tilde{B}_s$ is a $ (\shf_s)_{s\geq 0}$-Brownian motion, indeed 
\begin{align}
\sigma \int_0^{A^{-1}_s} a(X^{t,x}_r)^2 dr = \sigma  A(A^{-1}_s) =\sigma  s,
\end{align}
moreover setting $ A^{-1}_u = r$, 
\begin{align}
&\int_0^{A^{-1}_s} b(t - r, X^{t,x}_r) dr = \int_0^s  b(t - A^{-1}_u, X^{t,x}_{A^{-1}_u}) \frac{\partial A^{-1}_u}{\partial u} du\\
&\qquad =\int_0^s b(t - A^{-1}_u, Z^{t,x}_u) \frac{1}{\frac{\partial A}{\partial u}_{|_{A^{-1}_u}}} du.\label{eq_36b}
\end{align}
Note that $\frac{\partial A_u}{\partial u} = a(X^{t,x}_u)^2$, so 
\begin{equation}\label{eq_36c}
\frac{1}{\frac{\partial A}{\partial u}_{|_{A^{-1}_u}}} = \frac{1}{a(X^{t,x}_{A^{-1}_u})^2} =\frac{1}{a(Z^{t,x}_u)^2}.
\end{equation}
Finally, by \eqref{eq_36a}, \eqref{eq_36b}, \eqref{eq_36c},
\begin{align}
Z^{t,x}_s = x + \sigma \tilde{B}_s + \int_0^s b(t - A^{-1}_u, Z^{t,x}_u) \frac{1}{a(Z^{t,x}_u)^2}du +L^{t,x,Z}_s-U^{t,x,Z}_s,
\end{align}
for $0\leq s \leq A_{t}$. 
\end{proof}
Observe that the paths of the Brownian motion $\tilde{B}_s$ depend on $t,x$. Now we need some properties of the functions $A$ and $A^{-1}$.

\begin{rem}
Another possible method to reduce $X_s^{t,x}$ to a process with constant diffusion term would be to use the transform $R$ defined in \eqref{Rdef} which yields the process  $R(X_s^{t,x})$ in  \eqref{not:star1}. Then proceed as in Proposition \ref{Rem:3.4} to bound $R(X_s^{t,x})$ from above and below with a arithmetic Brownian motion. The drawback of this approach is that the new drift, instead of being just a scaled version of the original one, see \eqref{Z_def}, would also become distorted by the translation $-\frac{\s}{2}a'(x)$.
\end{rem}

\begin{cor}\label{Cor:1.42}
Assume Hypothesis \ref{hyp}. For any $t \geq 0$ and $ x\in \bar \shi$ the following holds:
\begin{enumerate}[(1)]
\item $ (c^d_a)^2 s\leq A^{t,x}_s\leq  (c_a^u)^2 s$, for $s\in [0,t]$;
\item $ \frac{u}{(c_a^u)^2}\leq A^{-1,t,x}_u \leq \frac{u}{(c^d_a)^2}$, for $u\in [0, A^{t,x}_{t}]$,
\end{enumerate}
where $c^d_a, c_a^u$ are the constants defined in Item 4) of Hypothesis \ref{hyp}.
\end{cor}
\begin{proof}
\begin{enumerate}[(1)]
\item The proof directly follows from the definition of $A_s$.
\item $A_s\leq  (c_a^u)^2 s$, so setting $u=A_s$, it follows that $u\leq (c_a^u)^2 s$. Then $s \geq \frac{u}{(c_a^u)^2}$ and recalling that $s= A^{-1}(A_s) = A^{-1}_u$,  the inequality of the left hand side is proved. The other inequality follows in the same way from $A_s \geq (c^d_a)^2 s $.
\end{enumerate}
\end{proof}

We now define some stopping times that will be used in the sequel. 
\begin{align}\label{hit_time}
&\tau^{t,x,X}_\partial:= \inf \{s \leq t :\ X^{t,x}_s \in \partial I \},\\
&\tau^{t,x,Z}_\partial:= \inf \{s\leq  A^{t,x}_{t}:\ Z^{t,x}_s \in \partial I \} ,\\
&\tau^{t,x,Z}_0:= \inf \{s\leq  A^{t,x}_{t}:\ Z^{t,x}_s =0 \},
&&\tau^{t,x,Z}_\zeta:= \inf \{s\leq  A^{t,x}_{t}:\ Z^{t,x}_s  = \zeta \},\\
&\tau^{\mu,x,S}_0:= \inf \{s:\ S^{\mu,x}_s =0 \},
&&\tau^{\mu,x,S}_\zeta:= \inf \{s:\ S^{\mu,x}_s  = \zeta \}.
\end{align}

The next Corollary allows us to compare the hitting times of the boundary $\partial I$ of the original process $X^{t,x}_s$ with the ones of the time changed process  $Z_s$ whose law is sometimes easier to estimate. Recall that $\partial \shi:= \{ 0\} \cup \{ \zeta \}$.

\begin{cor}\label{Cor:1.43}
Assume Hypothesis \ref{hyp}.  For any $t \geq 0$ and $ x\in \bar \shi$ the following holds:
\begin{enumerate}[1)]
\item $ (c_a^d)^2 \tau^{t,x,X}_\partial \leq \tau^{t,x,Z}_\partial $ a.s.;
\item $\P(\tau^{t,x,X}_\partial   \geq  \frac{s}{(c_a^d)^2})  \leq \P(\tau^{t,x,Z}_\partial  \geq s)$, for $s \in [0,(c_a^d)^2 t]$.
\end{enumerate}
\end{cor}

\begin{proof}
Fix $\omega \in \Omega$ such that $\tau^{t,x,Z}_\partial(\omega) <A^{t,x}_{t}(\omega)$ and let $s = s(\omega) \in [0, A^{t,x}_{t}[$ be such that $ Z^{t,x}_s(\omega) \in \partial I $ for the first time. By the definition of $Z$ in Lemma \ref{Lem:1.41}, this corresponds to $ X^{t,x}_{A^{-1}_s}(\omega) \in \partial \shi $. By Item 2) of Corollary \ref{Cor:1.42}, $A^{-1}_s(\omega) \leq \frac{s(\omega)}{(c_a^d)^2}$, so at time $\frac{s(\omega)}{(c_a^d)^2}$ the process  $X^{t,x}_{s}(\omega)$ has already hit $\partial I$ at least once. 
In other words, $\tau^{t,x,X}_\partial(\omega) \leq \frac{s(\omega)}{(c_a^d)^2}$ and $\tau^{t,x,Z}_\partial(\omega) = s(\omega)$.

If $\tau^{t,x,Z}_\partial(\omega) = A^{t,x}_{t}(\omega)$, that is, the process $Z$ does not hit the boundary, then the claim is also proved since also $X$ does not hit the boundary until time $t$ and, by Item 1) of Corollary \ref{Cor:1.42}, $  (c^d_a)^2 t \leq  A^{t,x}_{t}$ a.s.

Item 2) follows then from Item 1).
\end{proof}

The next Proposition puts all the above results together.
Fix $t\geq 0$ and define the following two drifts;
\begin{equation}\label{drift_def}
\mu_l := \inf_{s\in [0,t],x\in \bar\shi} \frac{b(t-s,x)}{a(x)^2},\qquad \mu_u := \sup_{s\in [0,t],x\in \bar\shi} \frac{b(t-s,x)}{a(x)^2}.
\end{equation}

\begin{prop}\label{Rem:3.4}
Assume Hypothesis \ref{hyp}. Let $t\geq 0$, $s\in [0,t]$ and $x\in \bar \shi$.
Then the following holds
\begin{align}
&\P(\tau^{t,x,X}_\partial \geq s)  \leq \P(\tau^{\mu_u,x,S}_0 \geq (c^d_a)^2 s),\label{bound_tauXu}\\
&\P(\tau^{t,x,X}_\partial \geq s)  \leq \P(\tau^{\mu_l,x,S}_\zeta \geq (c^d_a)^2 s). \label{bound_tauXd}
\end{align}
\end{prop}
\begin{proof}

The probability of hitting $\partial \bar \shi$ before a given time is higher than the probability of hitting one of the two boundaries before the same given time, hence, for  $s\in [0,A^{t,x}_{t}]$,
\begin{align}
&\P(\tau^{t,x,Z}_\partial  \geq s) \leq \P(\tau^{t,x,Z}_0   \geq s),   \label{bound_tauZ1}\\
&\P(\tau^{t,x,Z}_\partial  \geq s) \leq \P(\tau^{t,x,Z}_\zeta  \geq s) .\label{bound_tauZ1a}
\end{align}
The above holds, in particular, for $s\in [0, (c^d_a)^2 t]$ since $c^d_a t \leq A^{t,x}_{t}$ by Corollary \ref{Cor:1.42}.

By Proposition \ref{comp1} and the strong Markov property, the following holds
 \begin{align}
& Z^{t,x}_s \leq S^{\mu_u,x}_s, \, a.s. \qquad s\in [0, \tau^{t,x,Z}_0]\\
& Z^{t,x}_s \geq S^{\mu_l,x}_s, \, a.s.  \qquad s\in [0,\tau^{t,x,Z}_\zeta ],
\end{align}
where the process $S^{\mu_u,x}$ (resp $S^{\mu_l,x}$) is reflected at $\zeta$ (resp $0$) and defined in \eqref{eq:5.1a} (resp. \eqref{eq:5.1b}).
In particular, for $s\in [0,(c^d_a)^2 t] $,
\begin{align}
&\P(\tau^{t,x,Z}_0  \geq s)  \leq \P(\tau^{\mu_u,x,S}_0 \geq s), \label{bound_tauZ2}\\
&\P(\tau^{t,x,Z}_\zeta  \geq s)  \leq \P(\tau^{\mu_l,x,S}_\zeta \geq s).\label{bound_tauZ3}
\end{align}

In view of Item 2) of  Corollary \ref{Cor:1.43}, \eqref{bound_tauZ1} (resp. \eqref{bound_tauZ1a}),  and \eqref{bound_tauZ2} (resp. \eqref{bound_tauZ3}), for $s \in [0,(c^d_a)^2 t]$
\begin{align}
&\P(\tau^{t,x,X}_\partial \geq  \frac{s}{(c^d_a)^2})  \leq \P(\tau^{t,x,Z}_0 \geq s) \leq \P(\tau^{\mu_u,x,S}_0 \geq s),\\
&\P(\tau^{t,x,X}_\partial \geq  \frac{s}{(c^d_a)^2})  \leq \P(\tau^{t,x,Z}_\zeta \geq s) \leq \P(\tau^{\mu_l,x,S}_\zeta \geq s).
\end{align}
 By doing the change of variable $u = \frac{s}{(c^d_a)^2}$,  \eqref{bound_tauXu} and \eqref{bound_tauXd} follow.
\end{proof}

By \eqref{tauBound},  we have  a spectral  expansion for $\P(\tau^{\mu_u,x,S}_0 \geq s)$ and $\P(\tau^{\mu_l,x,S}_\zeta \geq s)$ which will be a key component in the proof of the main theorems.

\section{Differentiability of $X_s^{t,x}$.}\label{Sec:4}

In this Section we prove that $x \ra \frac{\partial X^{t,x}_s }{\partial x}  $ is continuous on $ \shi$.
A multidimensional statement appears in \cite{DeuZa}, Theorem~1, under the stronger assumptions that the drift is $C^1$ and time-homogeneous. Here we work in one dimension, allow time dependence in the drift, and assume that its spatial derivative has  finitely many discontinuities.

The argument of \cite{DeuZa} applies when the diffusion coefficient is constant, which is not our case. A direct time change does not resolve this either, because the Brownian motion $\tilde B$ constructed in Lemma~\ref{Lem:1.41} depends on $(t,x)$. Instead, we compose $X$ with the transformation $R$ defined below; so that the process $R(X)$ then has constant diffusion. Moreover, while \cite{DeuZa} considers reflection at a single boundary ($x=0$), our setting involves reflection at both endpoints; see \eqref{FC}.


Theorem~\ref{ContBnd} provides an explicit formula for $\frac{\partial}{\partial x} R\!\left(X^{t,x}_s\right)$ and shows that, after the first hitting time of the boundary by $X$ (denoted $\tau^{t,x,X}_\partial$), one has
$
\frac{\partial X^{t,x}_s}{\partial x}
=\frac{\partial R(X^{t,x}_s)}{\partial x} =0$
for all  $s\ge \tau^{t,x,X}_\partial$.
Throughout, we restrict to one-dimensional processes.


We define the function $R$ by
\begin{align}\label{Rdef}
R(x):= \frac{1}{\s} \int_0^x \frac{1}{a(z)}dz , \qquad x\in \bar \shi.
\end{align}
By It\^{o}' s formula, the process $R(X^{t,x}_s)$ satisfies, for $s\in [0,t]$,
\begin{align}
R(X^{t,x}_s) &= R(x)+ \int_0^s \left( \frac{b(t-r,X^{t,x}_r)}{\sigma a(X^{t,x}_r)} -\frac{\s}{2} a'(X^{t,x}_r)   \right)dr +B_s \label{eq:6.8R}\\
&\qquad +L^{R,t,x}_s-U^{R,t,x}_s,\\ 
R(X^{t,x}_s)&\in [0,R(\zeta)]\, a.s.\label{eq:6.9R}
\end{align}
Observe that by It\^{o}'s formula, $ L^{R,t,x}_s= \frac{1}{\s a(0)} L^{t,x}_s$ and  $ U^{R,t,x}_s= \frac{1}{\s a(\zeta)} U^{t,x}_s$.
Define $\tilde{a}(x):= a(R^{-1}(x)) $ and $\tilde{b}(x) := b(R^{-1}(x))$.  Equivalently, \eqref{eq:6.8R} can be rewritten as
\begin{align}\label{eq:6.8RR}
R(X^{t,x}_s) &= R(x)+ \int_0^s \left( \frac{\tilde{b}(t-r,R(X^{t,x}_r))}{\sigma \tilde{a}(R(X^{t,x}_r))}  -\frac{\s}{2} \tilde{a}'(R(X^{t,x}_r))   \right)dr +B_s\\
&\qquad  +L^{R,t,x}_s-U^{R,t,x}_s.
\end{align}
Define 
\begin{align}\label{not:star}
&\tilde{D}^{t}(s,x) := D^{t}(s,R^{-1}(x)),\qquad  D^{t}(s,x):= \frac{b(t-s,x)}{\sigma a(x)}  -\frac{\s}{2} a'(x),
\end{align}
then \eqref{eq:6.8RR} simplifies to
\begin{align}\label{not:star1}
R(X^{t,x}_s) = R(x)+ \int_0^s \tilde{D}^t(r,R(X^{t,x}_r)) dr +B_s  +L^{R,t,x}_s-U^{R,t,x}_s.
\end{align}
On the top of the hitting times defined in \eqref{hit_time}, we add those of $R(X^{t,x}_s)$. Let $R(\partial \shi) := \{0 \} \cup \{R(\zeta) \}$ (observe that $R(0) = 0$), and 
\begin{align}
\tau^{t,x,R}_\partial:= \inf \left\{s:\, R(X^{t,x}_s) \in  R(\partial \shi) \right\}.
\end{align}
Similarly to \eqref{hit_time}, define
\begin{align}
&\tau_0^{t,x,R}:= \inf \left\{s\leq t:\,  R(X^{t,x}_s) = 0  \right\} ,\  \tau_\zeta^{t,x,R}:= \inf \left\{s\leq t:\,  R(X^{t,x}_s) = R(\zeta)  \right\},\\
&\tau_0^{t,x,X}:= \inf \left\{s\leq t:\,  X^{t,x}_s = 0  \right\} ,\ \quad \tau_\zeta^{t,x,X}:= \inf \left\{s\leq t:\,  X^{t,x}_s = \zeta \right\}.
\end{align}
We remark that the result below holds for $x\in \shi$. The argument does not apply for $x\in \partial \shi$.

\begin{thm}\label{ContBnd}
For any $t \geq 0$ and $ x\in  \shi$ the following holds.
\begin{enumerate}[(1)]
 \item $\tau^{t,x,R}_\partial = \tau^{t,x,X}_\partial$ a.s.
In particular, a.s.  
\begin{equation}
\tau_0^{t,x,R}= \tau_0^{t,x,X},\qquad  \tau_\zeta^{t,x,R}= \tau_\zeta^{t,x,X}.  
\end{equation}

 \item
\begin{equation}
\frac{\partial R(X^{t,x}_s) }{\partial x} = R'(x)+ \int_0^s \frac{\partial R(X^{t,x}_r) }{\partial x}  \frac{\partial }{\partial x} \tilde{D}^{t}(r, \cdot)_{|R(X^{t,x}_r)} dr,
\end{equation}
a.s. for $s \in [0,\tau^{t,x,X}_\partial[ $; moreover
\begin{equation}
\frac{\partial R(X^{t,x}_s) }{\partial x} = \frac{\partial X^{t,x}_s }{\partial x} =0,
\end{equation}
a.s. for $s\in [ \tau^{t,x,X}_\partial , t]$, 
\item The functions
\begin{equation}
x \ra \frac{\partial R(X^{t,x}_s) }{\partial x},\qquad x \ra \frac{\partial X^{t,x}_s }{\partial x} 
\end{equation}
are a.s. continuous for  $s \in [0,t]$.

\item There exists two constants $C_1,C_2$ such that a.s. for  $s \in [0,t]$,
\begin{equation}
\frac{\partial R(X^{t,x}_s) }{\partial x} \leq C_1 \qquad \text{ and } \qquad  \frac{\partial X^{t,x}_s }{\partial x} \leq C_2.
\end{equation}
\end{enumerate}
\end{thm}

\begin{proof}

\textit{Proof of Item (1)}. This follows since $x \ra R(x)$ is bijective for $x \in [0,\zeta]$.

\textit{Proof of Item (2) and (3)}.
The proof of these items is divided in 5 steps. Steps from 1 to 4 follow closely the proof of \cite{DeuZa}, after taking into account the double barrier. The proof of Step 5 is different from \cite{DeuZa} and takes into account that the drift is not $C^{0,1}$ everywhere and has time dependence. 

\smallskip

\textit{Step 1}
We shall prove that there exists $\kappa >0$ such that for any $\ve >0$
 \begin{equation}\label{Gronw2}
(R(X^{t,x+\ve}_s)-R(X^{t,x}_s))^2\leq (R(x +\ve)-R(x))^2 \exp\left(2\kappa s\right), \, a.s., \, s\in [0,t].
\end{equation}
We follow the construction given in \textit{Step 1} of the proof of Theorem 1 of \cite{DeuZa} with some minor modifications.
Recalling  \eqref{not:star1},
for any $s \in [0,t] $, and every $\ve >0$, by It\^{o}'s formula 
\begin{align}
&d(R(X^{t,x+\ve}_s)-R(X^{t,x}_s))^2 \\
&= 2 (R(X^{t,x+\ve}_s)-R(X^{t,x}_s))\cdot( \tilde{D}^{t} (s,R(X^{t,x+\ve}_s)) - \tilde{D}^{t} (s,R(X^{t,x}_s)) )ds\\
&\qquad  +2 (R(X^{t,x+\ve}_s)-R(X^{t,x}_s))( dL^{R,t,x+\ve}_s-  dL^{R,t,x}_s ) \\
&\qquad -2(R(X^{t,x+\ve}_s)-R(X^{t,x}_s)) (dU^{R,t,x+\ve}_s - dU^{R,t,x}_s).
\end{align}
Observe that by Proposition \ref{comp} and the monotonicity of $x \ra R(x)$ entails that $R(X^{t,x+\ve}_s)\geq R(X^{t,x}_s)$, so
$ dL^{R,t,x+\ve}_s$ is positive if and only if $R(X^{t,x+\ve}_s)= R(X^{t,x}_s)=R(0)=0$. Similarly $ dU^{R,t,x}_s$ is positive if and only if $R(X^{t,x+\ve}_s)= R(X^{t,x}_s)=R(\zeta)$. Moreover,
\begin{align}
&-(R(X^{t,x+\ve}_s)-R(X^{t,x}_s))  dL^{R,t,x}_s  \leq 0\\
 &-(R(X^{t,x+\ve}_s)-R(X^{t,x}_s)) dU^{R,t,x+\ve}_s \leq 0,
\end{align}
 therefore
\begin{align}
&d(R(X^{t,x+\ve}_s)-R(X^{t,x}_s))^2 \\
&\quad \leq 2 (R(X^{t,x+\ve}_s)-R(X^{t,x}_s))\cdot( \tilde{D}^{t} (s,R(X^{t,x+\ve}_s)) - \tilde{D}^{t} (s,R(X^{t,x}_s)) )ds.
\end{align}
 By the Lipschitz property of $x\ra \tilde{D}^{t}(s,x)$, with Lipschitz constant uniform in $s$, there exists $\kappa$ such that 
\begin{align}
d(R(X^{t,x+\ve}_s)-R(X^{t,x}_s))^2 &\leq 2 \kappa (R(X^{t,x+\ve}_s)-R(X^{t,x}_s))^2.
\end{align}
An application of Gronwall's Lemma (see for example \cite[Chapter VII, Section X]{walter1998ode}) yields \eqref{Gronw2}.

Before \textit{Step 2.}, we introduce some notations that will be used in the rest of the proof.  Define $C := C_0 \cup C_\zeta \cup \{ t\}$, where $C_0 := \{s \in [0,t]:\, X^{t,x}_s=0 \}$ and $C_\zeta := \{s \in [0,t]:\, X^{t,x}_s=\zeta \}$. We also use the convention that $\inf \emptyset =  \{ t \}$.  Observe that, by construction, $$\inf C  =   \tau^{t,x,X}_\partial .$$
Without loss of generality, fix $\omega$ such that $\inf C = \inf C_0 <\inf C_\zeta$ and $\inf C < t $ (the other case $\inf C = \inf C_\zeta <\inf C_0$, and  $\inf C = t$ are treated in the same way), which corresponds to $\tau_0^{t,x,X}(\omega) < \tau_\zeta^{t,x,X}(\omega)$.  Choose a random $\beta >0$ such that $\tau_0^{t,x,X} < \tau_\zeta^{t,x,X}- \beta$. 

Let $(A_n)_{n\in \N}$ be the countable collection of  open connected components of the set $[0,t] \setminus C$. Since $A_n$ is open in $[0,t] $, fix $q_n \in A_n\cap \Q$. Finally, define $\s_s := \sup(C \cap [0,s])$, then $s\mapsto \s_s$ is locally constant on $[0,t] \setminus C$.
Observe that $\s_{q_n} = \inf A_n$, so $\s_{q_n}$ does not depend on the $q_n\in A_n$.

\smallskip

\textit{Step 2.}
Now we are going to show that for any $n \in \N$,
for all $s \in A_n$, $s \leq \tau_\zeta^{t,x,X}- \beta$ the following holds a.s.
\begin{align}\label{26}
L^{R,t,x}_s = L^{R,t,x}_{\s_{q_n}} =  [-R(x) -W_{\s_{q_n}}^{t,x}  ]^+.
\end{align}
 Set 
\begin{align}
W_s^{t,x}:= \int_0^s \tilde{D}^{t} (r,R(X^{t,x}_r)) dr + B_s.
\end{align}
Then, as in Equation (26) of Theorem 1 in \cite{DeuZa}, using Skorohod's Lemma (see \cite{ry3}, Chapter VI, Section 2), 
\begin{align}
L^{R,t,x}_s &= \sup_{r\leq s} [R(x) + W_r^{t,x} ]^-= \sup_{r\leq s} [(-R(x) - W_r^{t,x})\vee 0 ]\\
&= [\sup_{r\leq s} (-R(x) - W_r^{t,x}) ]\vee 0= [-R(x) +\sup_{r\leq s} ( - W_r^{t,x}) ]\vee 0\\
&=[-R(x) -\inf_{r\leq s}  W_r^{t,x} ]\vee 0 = [-R(x) -\inf_{r\leq s}  W_r^{t,x} ]^+.\label{26a}
\end{align}
By Girsanov's Theorem, there exists a probability measure $\tilde{\P}\ll \P$ under which $W^{t,x}_s$ is a Brownian motion. Since a Brownian Motion attains its minimum at a unique point over any interval of the form $[0,q]$ (see for example \cite{ry3}, Chapter III, Ex. (3.26) or \cite{kim1990}, Lemma 2.6), and the two probability measures are absolutely continuous, then also  $W^{t,x}_r$ attains its minimum at a unique point over $[0,q]$. Therefore the following random variable is well defined: 
$\a_{q}:=\arg\inf_{\{r\in [0,q]\}} \{ W^{t,x}_r\}$.
Using equation \eqref{26a}, for $q_n \in A_n$, $q_n \leq \tau_\zeta^{t,x,X}- \beta$,
\begin{align}
L^{R,t,x}_{q_n} = [-R(x) -  W_{\a_{q_n}}^{t,x} ]^+.
\end{align} 
\textit{Claim}.  For all $n$, $\s_{q_n}= \a_{q_n}$ a.s.

\textit{Proof of the claim}.
Fix $n$ and recall that $\sigma_{q_n} = \inf A_n$ and $\sigma_{q_n}\in C$. We divide the proof in three cases.

\textit{Case 1) Points of $C$ before $\sigma_{q_n}$.}
 For any $ c \in C$ with $c<\sigma_{q_n} $,  $R(X_{c}^{t,x}) = R(X_{\sigma_{q_n} }^{t,x})=0$ but $L_{c}^{R,t,x} < L_{\sigma_{q_n} }^{R,t,x}$ (since both $c$ and $\sigma_{q_n}$ belong to the support of the local time, which is increasing in these points). This results in  $W_{c}^{t,x} > W_{\sigma_{q_n} }^{t,x}$. 

\textit{Case 2) Points outside $C$ before $\sigma_{q_n}$.}
Also, for  any other $s\in [0,\sigma_{q_n} [ \setminus{C}$,  $W_{s}^{t,x} > W_{\sigma_{q_n} }^{t,x}$. 

\textit{Case 3) Points in $A_n$.}
For any $s\in A_n$, $W_{s}^{t,x} > W_{\sigma_{q_n} }^{t,x}$ because $R(X_{s}^{t,x}) > R(X_{\sigma_{q_n} }^{t,x})$ and $L_{s}^{R,t,x} = L_{\sigma_{q_n} }^{R,t,x}$ since $s$ and $q_n$ belong to the same excursion of the Brownian motion.

The three cases above prove that $\sigma_{q_n} = \alpha_n$ a.s.  Since the local time is constant on $A_n$, it also holds that $L_s^{R,t,x} = L_{q_n}^{r,t,x}$, for $s\in A_n$. Hence \eqref{26} is verified and  \textit{Step 2.} is proved.

\smallskip

\textit{Step 3.}
We are going to prove that, for a small enough perturbation of the initial value, $x+ \ve$, 
\begin{align}\label{27}
& L^{R,t,x+\ve}_{q_n} = [-R(x +\ve) -W^{t,x +\ve}_{\s_{q_n}}]^+\\
&\,= \left[ -R(x +\ve) - W^{t,x}_{\s_{q_n}} - \int_0^{\s_{q_n}} \tilde{D}^t(r, R(X^{t,x+\ve}_r)) - \tilde{D}^t(r,R(X^{t,x}_r) dr \right]^+,
\end{align}
where $\sigma_{q_n}$ is the same of \eqref{26} and, in particular, is independent of $\ve$. This is the key argument  to prove that,  after the first hitting time of the boundary, the derivative vanishes, see \eqref{28}.

The argument follows  closely the proof of Theorem 1, \textit{Step 4}, of \cite{DeuZa}.

Recalling that we are on the event $ \{\tau_0^{t,x,X} < \tau_\zeta^{t,x,X}\} $ and $s\ra X^{t,x}_s$ is continuous a.s.,  it is therefore possible to choose a  random $\Delta>0$ small enough such that for every $\ve \leq \Delta$, for all $s \in [	0,\tau_\zeta^{t,x,X}- \beta[$, both $U^{R,t,x + \ve}_{s} = 0$ and $U^{R,t,x}_{s}= 0$.

Since the law of $W_s^{t,x}$ over $s\in [0,t]$ is absolutely continuous w.r.t. the law of a Brownian motion, by Lemma 1 of \cite{DeuZa}, for every $q\in [0,t]  \cap \Q$ there exists a random variable $\gamma_q >0$ such that every small $\gamma_q$-Lipschitz perturbation of $W_s^{t,x}$ attains its minimum over $[0,q]$ only at $\s_q$. Moreover, by \eqref{Gronw2} and the Lipschitz property of  $x\ra \tilde{D}^t(s,x)$, for $q_n\in A_n$ and $q_n\in [0,t]  \cap \Q$, possibly after choosing a random $\Delta_n$, for every $\ve \leq \Delta_n \w \Delta$,
\begin{align}\label{KeyLemma}
\sup_{s\in[0,q_n] }|\tilde{D}^{t}(s, R(X^{t,x+\ve}_s)) - \tilde{D}^{t}(s,R(X^{t,x}_s))|\leq \gamma_{q_n}, 
\end{align}
which entails that, setting $f(s):=W^{t,x+\ve}_s -W^{t,x}_s$, 
\begin{align}
| f(s_1)-f(s_2) |\leq \gamma_{q_n}|s_1-s_2|,\qquad s_1,s_2\in [0,q_n].
\end{align}
Under this construction, by Lemma 1 in \cite{DeuZa} it follows that also $W^{t,x +\ve}_s$ a.s. attains its minimum over $[0,q_n]$ only  at $\s_{q_n}$ (observe that $\s_{q_n}$ is independent of $\ve \leq \Delta_n \w \Delta$) and, 
in particular, using \eqref{26a}, \eqref{27} is proved.

\smallskip 

\textit{Step 4.}
Here we are going to prove that for $\ve$ small enough
\begin{align}\label{28}
R(X^{t,x+\ve}_{s} ) = R(X^{t,x}_{s}),\qquad a.s.\ \ s \in [\inf C , t]
\end{align} 
 and 
\begin{equation}\label{29}
\xi^{t,x,R}_s(\ve) = \frac{R(x +\ve)-R(x)}{\ve}  + \frac{1}{\ve} \int_0^s   \tilde{D}^{t}(r, R(X^{t,x+\ve}_r))  - \tilde{D}^{t}(r, R(X^{t,x}_r))dr, \quad  a.s. \ s \in [0, \inf C[.
\end{equation}

\emph{Case 1: proof of \eqref{28}.} 

Let $q_n \in ] \inf C, \tau_\zeta^{t,x,X} -\beta[ \; \cap \;  A_n$ . Then   $R(x) + W^{t,x}_{\s_{q_n}} <0$, where the inequality is strict, and possibly after choosing a smaller random $\Delta'_n$, by \eqref{Gronw2} and the Lipschitz property of $x \ra \tilde D^t(r,x)$, for every $\ve \leq \Delta'_n  \wedge \Delta $, we also have
\begin{equation}\label{27a}
R(x +\ve) + W^{t,x}_{\s_{q_n}} + \int_0^{\s_{q_n}} [\tilde{D}^t(r, R(X^{t,x+\ve}_r)) - \tilde{D}^t(r,R(X^{t,x}_r))] dr <0.
\end{equation}
Then, by \eqref{26} and \eqref{27}
\begin{equation}\label{27aa}
L^{R,t,x+\ve}_{q_n}- L^{R,t,x}_{q_n} = - R(x +\ve) +R(x)  - \int_0^{\s_{q_n}} [ \tilde{D}^t(r, R(X^{t,x+\ve}_r)) - \tilde{D}^t(r,R(X^{t,x}_r)] dr,
\end{equation}
a.s. for $\ve < \Delta'_n \wedge \Delta_n \wedge \Delta$.
 Since $U_{q_n}^{R, t, x+ \ve} = U_{q_n}^{R, t, x} = 0 $, 
then  a.s. 
by \eqref{27aa} and  \eqref{not:star1},
\begin{align}\label{27ab}
R(X^{t,x+\ve}_{q_n}) - R(X^{t,x}_{q_n}) = \int_{\s_{q_n}}^{q_n} \tilde{D}^{t}(r, R(X^{t,x+\ve}_r))  - \tilde{D}^{t}(r, R(X^{t,x}_r)) dr.
\end{align}
Additionally we observe that, 
\eqref{27ab} holds if we replace $q_n$ with  any $s\in \, ]\sigma_{q_n}, \tau_\zeta^{t,x,X} -\beta[ \,  \cap  \, A_n$.
In particular, for $s\ra \sigma_{q_n}$, 
\begin{equation}\label{27ac}
R(X^{t,x+\ve}_{\s_{q_n}}) = R(X^{t,x}_{\s_{q_n}}), \qquad a.s. \ \ \ve < \Delta'_n \wedge \Delta_n \wedge \Delta.
\end{equation} 
By the strong Markov property, $ R(X^{t,x+\ve}_{s}) = R(X^{t,x}_{s})$  a.s. for any $s > \sigma_{q_n}$ and $\ve < \Delta'_n \wedge \Delta_n \wedge \Delta.$
Since this  holds for any $n$, it is possible to  choose a subsequence (depending on $\omega$), $(q_n)_{n\in \N}$ such that $\s_{q_n} \searrow \inf C $,  $(\ve_n)_{n\in \N}$ such that $\ve_n < \Delta'_n \wedge \Delta_n \wedge \Delta$. 

Because $\inf C$ is finite a.s., it follows that $\lim_{n\ra \infty} \gamma_{q_n} >0 $ a.s. (since $\gamma_{q_n}>0$ a.s. if $q_n$ is finite by Lemma 1 of \cite{DeuZa})  which in turn implies $\lim_{n\ra \infty}  \Delta_n >0 $ a.s. Moreover, by construction, also    $\lim_{n\ra \infty}  \Delta'_n >0 $. 
Hence,  \eqref{28} is verified.



\smallskip

\emph{Case 2, proof of \eqref{29}}. Let now  $ s < \inf C$,  then $\s_{s} = s$ and  
\begin{equation}\label{30a}
R(x) + W^{t,x}_{\s_{s}} >0.
\end{equation}
This leads to two consequences. First, by \eqref{26}, $L_s^{R,t,x} =0$ a.s. for $s\in[0,\inf C[$. Second, since the inequality in \eqref{30a} is strict,
arguing as at the beginning of \textit{Case 1} and by possibly choosing a different $\Delta'_n$, the expression on the left hand side of \eqref{27a} is also strictly positive for $\ve_n < \Delta'_n \wedge \Delta_n \wedge \Delta$ which,  by \eqref{27}, yields to $L^{R,t,x+\ve}_{s}=0$ a.s.  for $s\in [0,\inf C[$. Equivalently, defining 
\begin{align}
\xi^{t,x,R}_s(\ve) := \frac{R(X^{t,x+\ve}_s)-R(X^{t,x}_s)}{\ve} , \qquad s\in [0,\inf C[,
\end{align} 
\eqref{29} holds. 
Observe  that the case $\inf C = \inf C_\zeta <\inf C_0$ is handled in the same way. 

\smallskip

\textit{Step 5}.
In this step we prove that as $\ve \ra 0$, equation \eqref{29} does converge to a limit, which is the $x$-derivative of the process $R(X^{t,x}_s)$ and is a continuous process. This is the novel part of the proof.
\smallskip

\textit{Existence of the limit to \eqref{29} as $\ve \ra 0$.} 

Using the fundamental theorem of calculus in the form
\begin{align}\label{FTC}
f(b)-f(a) = (b-a) \int_0^1 f'(a+ \lambda(b-a)) d\lambda,
\end{align}
applied to $f(\cdot) = \tilde{D}^{t}(s,\cdot) $ with $a= R(X^{t,x}_s)$, $b= R(X^{t,x+\ve}_s)$, by \eqref{29},  a.s. for $s< \inf C$:
\begin{equation}\label{Eta}
\xi^{t,x,R}_s(\ve) = \frac{R(x +\ve)-R(x)}{\ve} + \int_0^s  \xi^{t,x,R}_r(\ve)  \int_0^1 \frac{\partial }{\partial x} \tilde{D}^{t}(r, \cdot)_{| R(X^{t,x}_r) + \lambda (R(X^{t,x+\ve}_r)-R(X^{t,x}_r))} d\lambda dr ,
\end{equation}
and by \eqref{28},  for  $ s \geq \inf C $ (equivalently $s\geq \tau_\partial^{t,x,X}$), 
\begin{align}\label{Eta1}
\xi^{t,x,R}_s(\ve) &=  0.
\end{align}
By standard theory of ODE, the solution to \eqref{Eta} satisfies 
\begin{equation}\label{dRSol}
\xi^{t,x,R}_s(\ve) =  \frac{R(x +\ve)-R(x)}{\ve} \exp\left(\int_0^s \int_0^1 \frac{\partial }{\partial x}  \tilde{D}^{t}(r, \cdot)_{| R(X^{t,x}_r) + \lambda (R(X^{t,x+\ve}_r)-R(X^{t,x}_r))} d\lambda dr \right).
\end{equation}
Since $\frac{\partial }{\partial x} \tilde{D}^{t}(s,x)$ is bounded for any $x \in \bar \shi$, $s\in [0,t]$,  the Lebesgue Dominated Convergence Theorem applies to the right hand side of \eqref{dRSol}:
\begin{align}
&\lim_{\ve \ra 0 }\exp\left(\int_0^s \int_0^1 \frac{\partial }{\partial x}  \tilde{D}^{t}(r, \cdot)_{| R(X^{t,x}_r) + \lambda (R(X^{t,x+\ve}_r)-R(X^{t,x}_r))} 1_{\{X_r^{t,x} \notin \shg \}} d\lambda dr \right)\\
&\qquad =  \exp\left( \int_0^s  \frac{\partial }{\partial x} \tilde{D}^{t}(r, \cdot)_{|R(X^{t,x}_r)} 1_{\{X_r^{t,x} \notin \shg \}} dr \right) =  \exp\left( \int_0^s  \frac{\partial }{\partial x} \tilde{D}^{t}(r, \cdot)_{|R(X^{t,x}_r)} dr \right),
\end{align}
where the last equality follows because  $\shg_R := \cup_{i=1}^g R(\hat x_i)$, the set of discontinuities of $x \mapsto \frac{\partial }{\partial x} \tilde{D}^{t}(s,x)$, (recall that, by Hypothesis \ref{hyp}, $ x\mapsto \frac{\partial }{\partial x}  D^t(s,x)$ is discontinuous on $\shg$) has zero Lebesgue measure.

Therefore, since the limit on the right hand side of \eqref{dRSol} exists, also the limit on left hand side exists and is equal to $ \frac{\partial R(X^{t,x}_s) }{\partial x}^+$, hence the following holds:
\begin{align}\label{1.108a}
 \frac{\partial R(X^{t,x}_s) }{\partial x}^+ = R'(x)+ \int_0^s \frac{\partial R(X^{t,x}_r) }{\partial x}^+  \frac{\partial }{\partial x} \tilde{D}^{t}(r, \cdot )_{|R(X^{t,x}_r)} dr,
\end{align}
a.s., $s\in[0, \inf C[$.
 Defining $\xi^{t,x,R}_s(-\ve) := \frac{R(X^{t,x}_s) -R(X^{t,x-\ve}_s)}{\ve}$ and repeating all the steps above, we see that the left derivative exists and coincides with the right one, so $ \frac{\partial R(X^{t,x}_s) }{\partial x}$ is well defined. After time $\tau^{t,x,X}_\partial $, the derivative is null by \eqref{Eta1}.
 This proves Item (2).
 
\smallskip \textit{Continuity of the limit. }

We first prove continuity of $x \ra \frac{\partial R(X^{t,x}_s) }{\partial x}$, for $s\in [0,\tau^{t,x,X}_\partial[$.  We write, using \eqref{1.108a},
\begin{align}
 &\frac{\partial R(X^{t,x+\ve}_s) }{\partial x}-   \frac{\partial R(X^{t,x}_s) }{\partial x} = R'(x+\ve) -R'(x)\\
&\qquad  +  \int_0^s \frac{\partial R(X^{t,x+\ve}_r) }{\partial x}  \frac{\partial }{\partial x} \tilde{D}^{t}(r, \cdot)_{| R(X^{t,x+\ve}_r )} dr- \int_0^s \frac{\partial R(X^{t,x}_r) }{\partial x}  \frac{\partial }{\partial x} \tilde{D}^{t}(r, \cdot)_{|R(X^{t,x}_r )} dr\\
&=  R'(x+\ve) -R'(x)  \\
&\qquad +  \int_0^s \left( \frac{\partial R(X^{t,x+\ve}_r) }{\partial x} -\frac{\partial R(X^{t,x}_r) }{\partial x}\right) \frac{\partial }{\partial x} \tilde{D}^{t}(r, \cdot)_{| R(X^{t,x+\ve}_r )} dr\\
&\qquad + \int_0^s \frac{\partial R(X^{t,x}_r) }{\partial x} \left( \frac{\partial }{\partial x} \tilde{D}^{t}(r, \cdot)_{|R(X^{t,x+\ve}_r)} -\frac{\partial }{\partial x} \tilde{D}^{t}(r, \cdot)_{| R(X^{t,x}_r )}\right) dr.
\end{align}

For sake of readability, let us define 
\begin{align}
\partial R_s:=  \frac{\partial R(X^{t,x}_s) }{\partial x}  , 
&& \Gamma_\ve \partial R_s&:= \frac{\partial R(X^{t,x+\ve}_s) }{\partial x}-   \frac{\partial R(X^{t,x}_s) }{\partial x}\\
\partial \tilde{D}^\ve_s:= \frac{\partial }{\partial x} \tilde{D}^{t}(s, \cdot )_{|R(X^{t,x+\ve}_s)},
&& \Gamma_\ve \partial \tilde{D}_s &:= \frac{\partial }{\partial x} \tilde{D}^{t}(s, \cdot)_{| R(X^{t,x+\ve}_s )} -\frac{\partial }{\partial x} \tilde{D}^{t}(s, \cdot)_{|R(X^{t,x}_s )},
\end{align}
 so that, by standard theory of ODE,
\begin{equation}\label{conv_drx}
\Gamma_\ve \partial R_s= \exp\left(\int_0^s \partial \tilde{D}^\ve_r dr\right)\left( R'(x+\ve)-R'(x) + \int_0^s  \exp\left(-\int_0^r  \partial \tilde{D}^\ve_u du\right) \cdot \partial R_r \cdot \Gamma_\ve \partial \tilde{D}_r dr \right). 
\end{equation}
Recalling \eqref{Gronw2}, define $K(\ve):=\sup_{x\in \shi} |R(x+\ve)-R(x)|\exp\left(\kappa t\right)$. By Hypothesis \ref{hyp}, the set of discontinuities of $\frac{\partial }{\partial x} \tilde{D}^{t}(s, x )$ is finite and equal to $\shg_R = \cup_{i=1}^{g} R(\hat x _i)$. Let $\hat x_0 =0$ and $\hat  x_{g+1}  = \zeta$, and  
\begin{align}\label{eq_3.24a}
& \shg^{\ve,c}_{1,R}:=  [R(\hat x_{0}), R(\hat x_{1}) - K(\ve)[; &&  \shg^{\ve,c}_{g+1,R}:=  ]R(\hat x_{g}) + K(\ve), R(\hat x_{g+1}) ], \\
&\shg^{\ve,c}_{i,R}:=  ]R(\hat x_{i-1}) + K(\ve), R(\hat x_{i}) - K(\ve)[, \ i \in \{2,\dots,g \},\\
& \shg^c_{1,R}:=  [R(\hat x_{0}) ,R( \hat x_{1}) [; &&  \shg^{c}_{g+1,R}:=  ]R(\hat x_{g}) , R( \hat x_{g+1}) ], \\
&\shg^c_{i,R} :=  ]R(\hat x_{i-1}) , R(\hat x_{i})[ ,  \ i \in \{2,\dots,g \} ,
\end{align}
 and their union sets $$ \shg^{\ve,c}_R:= \cup_{i=1}^{g+1} \shg^{\ve,c}_{i,R}, \qquad   \shg^c_R:= \cup_{i=1}^{g+1} \shg^c_{i,R}.$$
Now define $A_1(s) := \{\omega:  R(X^{t,x}_s) \in \shg^{\ve,c}_R \}$ and $A_2(s) := \{\omega: R(X^{t,x+\ve}_s) \in \shg^c_R  \}$. 
Then, for a small enough $\ve >0$,  $1_{A_2(s) \cap A_1(s)} = 1_{A_1(s)} $ a.s. for $s\in [0,t]$ because the base path and the perturbed path remain in the same open interval.
Hence  
\begin{equation}\label{eq_3.24b}
\lim_{ \ve\ra 0}\Gamma_{\ve} \partial \tilde{D}_s 1_{A_2(s) \cap A_1(s)}  = \lim_{ \ve\ra 0}\Gamma_{\ve} \partial \tilde{D}_s 1_{A_1(s)}  = 0 \qquad \quad\text{a.s. } s\in [0,\tau^{t,x,X}_\partial[.
\end{equation}
Therefore by \eqref{eq_3.24b}
\begin{equation}\label{eq_3.28}
\lim_{\ve \ra 0} \int_0^s  \exp\left(-\int_0^r  \partial \tilde{D}^\ve_u du\right) \cdot \partial R_r \cdot \Gamma_\ve \partial \tilde{D}_r \, 1_{A_1(r)} dr = 0 \quad \text{a.s.  } s\in [0,\tau^{t,x,X}_\partial[.
\end{equation}
where we could apply the dominated convergence theorem because  $\Gamma_\ve \partial \tilde{D}_r$  and $ \exp\left(-\int_0^r  \partial \tilde{D}^\ve_u du\right)$ are bounded by Hypothesis \ref{hyp} and $\partial R_r$ is bounded by \eqref{bound_drx}.

On the complement of $A_1$, $\{R(X^{t,x}_r)\in (\shg^{\varepsilon,c}_R)^c\}$, the limit of $\Gamma_\varepsilon \partial \tilde D_r$ may fail to exist.
Moreover, since $\lim_{\varepsilon\to 0} (\shg^{\varepsilon,c}_R)^c=\shg_R $ and the time spent by $R(X^{t,x}_\cdot)$ at the finite set $\shg_R$ has Lebesgue measure zero, another application of dominated convergence yields
 $\lim_{\ve \ra 0} \int_0^t  1_{\{\,R(X^{t,x}_r)\in (\shg^{\ve,c}_R)^c\,\}}\,dr = 0$ a.s. 
Thus
 \begin{equation}\label{eq_3.28a}
\lim_{\ve \ra 0}  \int_0^s  \exp\left(-\int_0^r  \partial \tilde{D}^\ve_u du\right) \cdot \partial R_r \cdot \Gamma_\ve \partial \tilde{D}_r  1_{\{\,R(X^{t,x}_r)\in (\shg^{\ve,c}_R)^c\,\}} dr = 0 \quad \text{a.s.  }  \, s\in [0,\tau^{t,x,X}_\partial[.
\end{equation}
Combining \eqref{conv_drx}, \eqref{eq_3.28}, and \eqref{eq_3.28a} we conclude that
$$
\Gamma_\varepsilon \partial R_s \longrightarrow 0
\quad\text{a.s.  } s\in [0,\tau^{t,x,X}_\partial[.
$$
 This proves the continuity of $x \ra \frac{\partial R(X^{t,x}_s) }{\partial x}$ for $s\in [0, \tau^{t,x,X}_\partial[$.
 By \eqref{28}, for small enough $\ve$,  $\tau^{t,x,X}_\partial = \tau^{t,x+\ve,X}_\partial$ a.s., hence  for $s\in [\tau^{t,x,X}_\partial,t]$, continuity holds as well because $\frac{\partial R(X^{t,x}_s) }{\partial x} =\frac{\partial R(X^{t,x+\ve}_s) }{\partial x}= 0$.
By the chain rule formula,
\begin{align}\label{ChainRule}
\frac{\partial R(X^{t,x}_s) }{\partial x}  = R'(X^{t,x}_s) \frac{\partial X^{t,x}_s}{\partial x},
\end{align}
and since $ x \ra R'(X^{t,x}_s)$ continuous, $x \ra \frac{\partial X^{t,x}_s }{\partial x}$ is  also continuous on $\bar \shi$. 
Item (3) is then proved.

\smallskip

\textit{Proof of Item (4)}.

It is enough to observe that, by Hypothesis \ref{hyp}, for any $(s,x) \in [0,t] \times \bar \shi$, there exists a constant $c >0$ such that  $\frac{\partial }{\partial x} \tilde{D}^{t}(s, x) \leq c $. Then, by  Gronwall's Lemma (see for example \cite[Chapter VII, Section X]{walter1998ode}), applied to \eqref{1.108a},
\begin{align}\label{bound_drx}
\left| \frac{\partial R(X^{t,x}_s) }{\partial x} \right| \leq R'(x) \exp\left(c s\right)
\end{align}
The bound on $\frac{\partial X^{t,x}_s }{\partial x}$ follows by the chain rule formula \eqref{ChainRule} observing that $R'(x) = \frac{1}{\s a(x)}$ is also bounded by Hypothesis \ref{hyp}.
	\end{proof}


\section{Explicit pointwise bound on $\left|\frac{\partial V}{\partial x}\right|$}\label{Sec:5}

In this section we prove the main result, Theorem \ref{final_thm}, which is a direct consequence of Theorems \ref{316_bound_a} and \ref{317_bound}. We also prove Corollary \ref{lim_infty} which provides a sufficient condition for $\left|\frac{\partial V}{\partial x}\right|$ to be bounded as $t\ra \infty$.

We first need a preliminary Lemma in which we provide a preliminary bound on the $x$-derivative of $V$, making use of the Feynman-Kac formula of Proposition \ref{FeyC} and  Item 2. of Theorem \ref{ContBnd}.
Define
\begin{equation}\label{cost_def}
c_D=\sup_{r \in [0,t],x \in \bar \shi} \left|\frac{\partial }{\partial x} D^t(r,x)\right|, \qquad c_{h} := \sup_{r \in [0,t],x \in \bar \shi} |h'(r,x)|, \qquad c_{g} := \sup_{x \in \bar \shi} |g'(x)|,
\end{equation}
and we recall that $c^d_a, c_a^u$ were  introduced in Hypothesis \ref{hyp}.

\begin{lem}\label{LThetaBound}
Assume Hypothesis \ref{hyp} is verified. Then for any $t\geq 0$ and $x\in \shi$,
\begin{align}
\left| \frac{\partial V}{\partial x}(t,x) \right| &\leq   \frac{c^u_a }{ c^d_a} \left( c_{h} \int_0^{t}  \exp\left(\sigma \, c_D \, c^u_a  \,   s\right) \P\left(s \leq \tau_\partial^{t,x,X}\right) ds \right.\\
&\qquad \left.+ c_{g} \,  \exp\left(\sigma \, c_D \, c^u_a  \,   t\right) \P(t \leq \tau_\partial^{t,x,X})\right) .
\end{align}
\end{lem}

\begin{proof}
Without losing the generality, we assume $h$ does not depend on $t$.
Using Proposition \ref{FeyC},
\begin{equation}\label{eq_510a1}
\frac{\partial V}{\partial x}(t,x) =  \lim_{\ve \ra 0} \E \left[ \int_0^{t}  \frac{h(X^{t,x + \ve}_s)- h(X^{t,x}_s) }{\ve} ds\right] + \E\left[\frac{ g(X^{t,x + \ve}_{t}) -g(X^{t,x}_{t})}{\ve}\right], 
\end{equation}
and by the fundamental theorem of calculus in the form \eqref{FTC}, for every $s$,   
\begin{equation}\label{FTCA1}
\frac{h(X^{t,x + \ve}_s)- h(X^{t,x}_s) }{\ve} = \frac{X^{t,x + \ve}_s - X^{t,x}_s}{\ve}  \int_0^1 h'(X^{t,x}_s+\lambda (X^{t,x+\ve }_s- X^{t,x}_s)) d\lambda.
\end{equation}
By Theorem \ref{ContBnd}, $\frac{X^{t,x + \ve}_s-X^{t,x}_s}{\ve} \ra \frac{\partial X^{t,x}_s}{\partial x}$ a.s. , $s\in [0,t]$, and by the continuity of $x \ra h'(x)$ (outside of the finite set $\shg$) and $ x \ra X^{t,x}_s$, 
\begin{equation}
\frac{X^{t,x + \ve}_s-X^{t,x}_s}{\ve} \int_0^1 h'(X^{t,x}_s+\lambda (X^{t,x+\ve }_s - X^{t,x}_s)) d\lambda \ra \frac{\partial X^{t,x}_s}{\partial x} h'(X^{t,x}_s) \qquad a.s.
\end{equation}
  Moreover, the right hand side of \eqref{FTCA1}  is bounded in $L^1([0,t]\times \Omega)$ therefore the Lebesgue Dominated Convergence Theorem applies. After using the same argument for the second term on the r.h.s of \eqref{eq_510a1}, we obtain 
\begin{equation}\label{5.212}
\frac{\partial V}{\partial x}(t,x)=  \E \left[ \int_0^{t}  h'(X^{t,x}_s) \frac{\partial X^{t,x}_s }{\partial x} ds \right] + \E\left[ g'(X^{t,x}_{t})\frac{\partial X^{t,x}_{t} }{\partial x}\right].
\end{equation}
By \eqref{ChainRule}, since $\frac{1}{R'(x)} \leq  \sigma \, c^u_a $, then $\left| \frac{\partial X^{t,x}_s }{\partial x}\right|  \leq \sigma \, c^u_a \, \left|\frac{\partial R(X^{t,x}_s)}{\partial x}\right|$. Moreover, by Item (2) of Theorem \ref{ContBnd} $\frac{\partial X^{t,x}_s }{\partial x} = 0$ for $s\geq \tau_\partial^{t,x,X}$, hence
\begin{align}\label{5.212a}
\left| \frac{\partial V}{\partial x}(t,x) \right|  &\leq  \sigma \, c^u_a \, c_{h} \, \E \left[ \int_0^{t} \left|  \frac{\partial R(X^{t,x}_s)}{\partial x} \right| 1_{\{s< \tau_\partial^{t,x,X}\}}  ds \right]\\
& \qquad +  \sigma \, c^u_a \, c_g \E \left[\left|  \frac{\partial R(X^{t,x}_t)}{\partial x} \right| 1_{\{t< \tau_\partial^{t,x,X}\}}\right].
\end{align}
We now prove a bound for $\left| \frac{\partial R(X^{t,x}_s)}{\partial x}\right|$ which does not depend on $\omega$. We will sometimes write 
$X_s$ in place of $X^{t,x}_s$ when this does not generate confusion. 
Theorem \ref{ContBnd} entails that, defining $Y_s:= \frac{\partial X_s}{\partial x}$, $Y^R_s:= \frac{\partial R(X_s)}{\partial x}$,  for $s < \tau^{t,x,X}_\partial$,
\begin{align}
Y^R_s &=  R'(x) + \int_0^s \frac{\partial }{\partial x} \tilde{D}^t(r,\cdot)_{|R(X_r)} Y^R_r dr\label{1.191}
\end{align}
where we recall that $D^t$ and $\tilde{D}^t$ are defined in \eqref{not:star1}.
By the chain rule formula, 
\begin{align}
\frac{\partial }{\partial x} \tilde{D}^t(s,x)|_{R(X_s)}= \frac{\partial }{\partial x} D^t(s,x)_{|R^{-1}(R(X_s)) } (R^{-1})'(R(X_s)). 
\end{align}
Observing that 
\begin{align}
(R^{-1})'(R(X_s))= \frac{1}{R'(R^{-1}(R(X_s)))} =  \frac{1}{R'(X_s)} = \sigma \, a(X_s),
\end{align}
equation \eqref{1.191} becomes
\begin{align}
& Y^R_s=R'(x) + \sigma \int_0^s \frac{\partial }{\partial x} D^t(r,\cdot)_{|X_r}  a(X_r) Y^R_r dr,
\end{align}
hence
\begin{align}
& \left|Y^R_s \right| \leq \left| R'(x) \right|+ \sigma \int_0^s  \left| \frac{\partial }{\partial x} D^t(r,\cdot)_{|X_r}  a(X_r) Y^R_r \right|  dr,
\end{align}
Since  $ \left| R'(x) \right| \leq \frac{1}{\sigma \, c^d_a}$, 
we can rewrite the above as
\begin{align}\label{Y_bound1}
 \left|Y^R_s \right|  \leq  \frac{1}{\sigma \, c^d_a} + c_D \, c^u_a \,   \sigma \int_0^s  \left|Y^R_r \right|  dr.
\end{align}
Gronwall's Lemma (see for example \cite[Chapter VII, Section X]{walter1998ode}) applied to \eqref{Y_bound1} yields
\begin{align}\label{Ybound}
& \left|Y^R_s \right|  \leq  \frac{1}{\sigma \, c^d_a}  \,  \exp\left( \sigma \, c_D \, c^u_a  \,   s\right), \qquad s <\tau^{t,x,X}_\partial.
\end{align}
Finally, by \eqref{5.212a} and \eqref{Ybound}  the conclusion follows.
\end{proof}


It is possible to bound $\P(s \leq \tau^{t,x,X}_\partial)$ in two ways, see \eqref{bound_tauXu}, \eqref{bound_tauXd}.  The bounds \eqref{dv_bound_1}, \eqref{eq_380a} and \eqref{eq_380} hold for $\mu_u \geq 0$ and use \eqref{bound_tauXu}, while  \eqref{dv_bound_2}  and  \eqref{eq_380aa_m} are valid in the case $\mu_u <0$ and make use of \eqref{bound_tauXd} and  the reflection argument of Remark \ref{mu_neg_rem}.

Proposition \ref{prop_g_bound} and Theorem \ref{316_bound_a} provide  estimates which are accurate near $x=0$, while  Proposition \ref{prop_d_bound} and Theorem \ref{317_bound} exhibit estimates which are accurate near $x=\zeta$.

We recall that $\mu_u$ and $\mu_l$ are defined in Proposition \ref{Rem:3.4}, while $\lambda_n^{u,\mu_u}$ (resp. $\lambda_n^{l,\mu_l}$) is introduced in \eqref{lambda_u_def} (resp. \eqref{lambda_l_def}). 


\begin{prop}\label{prop_g_bound}
Assume Hypothesis \ref{hyp} is verified. For any $t\geq 0$ and $x\in  \shi$, the following holds.
\begin{enumerate}[(1)]
\item If $\mu_u \geq 0$ then 
\begin{equation}\label{dv_bound_1}
\left| \frac{\partial V}{\partial x}(t,x) \right|  \leq  \frac{c^u_a }{ c^d_a }   \sum_{n=1}^\infty  \left| c^{u,\mu_u}_{n,0}(x) \right|  \shn( t,\lambda^{u,\mu_u}_n,\sigma), 
\end{equation}
\item otherwise, if $\mu_u < 0$, then 
\begin{equation}\label{dv_bound_2}
\left| \frac{\partial V}{\partial x}(t,x) \right|  \leq  \frac{c^u_a }{ c^d_a }  \sum_{n=1}^\infty  \left| c^{l,-\mu_u}_{n,\zeta}(\zeta- x) \right|\shn( t,\lambda^{l,-\mu_u}_n,\sigma), 
\end{equation}
\end{enumerate}
where
\begin{align}
&\shn( t,\lambda,\sigma):= c_h\shl(t,\lambda,\sigma)  + c_g \exp\left(\shd(\lambda,\sigma)t\right), \label{N_def_a}\\
&\shl(t,\lambda,\sigma) :=  \int_0^t  \exp\left( \shd(\lambda,\sigma) s\right)  \; ds,\label{L_def_a}\\
&\shd(\lambda,\sigma) := \sigma \, c_D \, c^u_a  -\lambda \, (c^d_a)^2, \label{D_def_a}
\end{align}
and  $\lambda^{u,\mu_u}_{n}$ (resp. $\lambda^{l,-\mu_u}_{n}$) is defined in \eqref{lambda_u_def} (resp. \eqref{lambda_l_def}) and $c_{n,0}^{u,\mu_u}$ (resp. $c^{l,-\mu_u}_{n,\zeta}$) are as in \eqref{cnDef1} (resp. \eqref{cnDef2}).
\end{prop}

\begin{proof}
We first prove  Item (1) with $g=0$ since the general case follows easily.
By Lemma \ref{LThetaBound} and \eqref{bound_tauXu},
\begin{equation}\label{eq342}
\left| \frac{\partial V}{\partial x}(t,x) \right| \leq\frac{c^u_a \, c_{h}}{ c^d_a}  \left|  \int_0^{t}  \exp\left(\sigma \, c_D \, c^u_a  \,   s\right) \P\left(\tau_0^{\mu_u,x,S} \geq (c^d_a)^2 s \right) ds \right|.
\end{equation}
We now use \eqref{tauBound} with the $c^{u,\mu_u}_{n,0}$ as in \eqref{cnDef1}, with $y=0$, to compute the spectral expansion of $\P\left(\tau_0^{\mu_u,x,S} \geq (c^d_a)^2 s \right)$.
By so doing, \eqref{eq342} becomes
\begin{equation}\label{eq343}
\left| \frac{\partial V}{\partial x}(t,x) \right| \leq  \frac{c^u_a \, c_{h}}{ c^d_a} \, \left| \int_0^{t} \left( \exp\left(\sigma \, c_D \, c^u_a  \,   s\right)  \sum_{n=1}^\infty  c_{n,0}^{u,\mu_u}(x) \exp\left(-\lambda^{u,\mu_u}_n \, (c^d_a)^2 s\right) \right)ds \right|.
\end{equation}

We now show that we can invert the sum with the integral.
For $n>1$, by \eqref{lambda_u_def}, $\lambda^{u,\mu_u}_n = \frac{(\alpha^{u,\mu_u}_n)^2\s^2}{2 \zeta^2} + \frac{\mu_u^2}{2\s^2}$.
If  $\frac{\s^2}{\mu_u \zeta} \geq 1$, for $n$ large enough, by Remark \ref{rem_34}, the following holds true
\begin{equation}\label{eq_514aa}
\left| c_{n,0}^{u,\mu_u}(x) \right|  \leq  4  \frac{ x}{\zeta} \exp\left(-\frac{\mu_u}{\sigma^2} x \right).
\end{equation}
Trivially, if  $\frac{\s^2}{\mu_u \zeta} < 1$ by  Item (2)  Proposition \ref{prop_cn_bound}, \eqref{eq_514aa} is also verified.

By \eqref{eq_514aa} and since $\alpha^{u,\mu_u}_n \in [(n-1)\pi, (n-1)\pi + \pim[$, see   \eqref{lambda_def}, for $n$ large enough, the following bound holds
\begin{equation}\label{eq_510a}
\left| c_{n,0}^{u,\mu_u}(x) \exp\left(-\lambda^{u,\mu_u}_n \, (c^d_a)^2 s\right) \right| \leq      4  \frac{ x}{\zeta} \exp\left(-\frac{\mu_u}{\sigma^2} x\right)  \exp\left(-\frac{ (n-1)^2\pi^2 \sigma^2 }{2\zeta^2} \, (c^d_a)^2 s\right)
\end{equation}
Therefore, the Lebesgue dominated convergence theorem applies and \eqref{eq343} becomes
\begin{equation}\label{eq343a}
\left| \frac{\partial V}{\partial x}(t,x) \right| \leq  \frac{c^u_a \, c_{h}}{ c^d_a} \, \left|  \sum_{n=1}^\infty \left( c_{n,0}^{u,\mu_u}(x)  \int_0^{t}\exp\left((\sigma \, c_D \, c^u_a  -\lambda^{u,\mu_u}_n \, (c^d_a)^2 )   s\right)   ds \right)\right|.
\end{equation}

Item (2). By the reflection principle of Remark \ref{mu_neg_rem},
$$\P\left(\tau_0^{\mu_u,x,S} \geq (c^d_a)^2 s \right) = \sum_{n=1}^\infty  c_{n,\zeta}^{l,-\mu_u}(\zeta-x) \exp\left(-\lambda^{l,-\mu_u}_n \, (c^d_a)^2 s\right) ,$$
then argue as in the proof of Item (1) above and use Item (3) of Proposition \ref{prop_cn_bound} to derive \eqref{eq_510a}.
\end{proof}

Now, we  derive a closed-form bound on $\sum_{n=1}^\infty   \shn( t,\lambda^{u,\mu_u}_n,\sigma) $ depending only on the coefficients and the first spectral eigenvalue, hence computable explicitly.

We recall that $\alpha_n^{u,\mu_u}$ (resp. $\alpha_n^{l,-\mu_u}$) is introduced in \eqref{a_def_eq} (resp. \eqref{a_def_eq_1}), moreover, $ \lambda_n^{u,\mu_u}$ (resp. $ \lambda_n^{l,-\mu_u}$) is introduced in \eqref{lambda_u_def} (resp. \eqref{lambda_l_def}).

\begin{lem}\label{lem_bound_1}
Assume Hypothesis \ref{hyp} is verified. 
Let 
 $(\lambda_n)_{n\in \N}$ be either 
\begin{align}
&\text{\emph{Case A:}} \quad (\lambda^{u,\mu_u}_{n})_{n\in \N},  \quad  \lambda^{u,\mu_u}_{n} = \frac{\mu_u^2}{2\sigma^2} \;+\; \frac{\sigma^2}{2\zeta^2}\,(\alpha_n^{u,\mu_u})^2, \quad &\text{if} \ \mu_u \geq 0,\\ 
&\text{\emph{Case B:}}\quad  (\lambda^{l,-\mu_u}_{n})_{n\in \N}, \quad \lambda^{l,-\mu_u}_{n} = \frac{\mu_u^2}{2\sigma^2} \;+\; \frac{\sigma^2}{2\zeta^2}\,(\alpha_n^{l,-\mu_u})^2, \quad &\text{if} \ \mu_u < 0;
\end{align}
and  $\alpha_n = \alpha_n^{u,\mu_u}$ (resp. $\alpha_ n = \alpha_n^{l,-\mu_u}$) if \emph{Case A} (resp. \emph{Case B}) holds.
Write $$ \shd(\lambda_n,\sigma)= \sigma \, c_D \, c^u_a  -\lambda_n \, (c^d_a)^2=c_{\mu_u}-k\alpha_n^2,$$ with
\begin{equation}\label{cmuu}
c_{\mu_u} := \sigma c_D c^u_a - \frac{(c^d_a)^2\mu_u^2}{2\sigma^2},
\qquad
k := \frac{(c^d_a)^2\sigma^2}{2\zeta^2}.
\end{equation}

Then, the following bounds hold for all $t \ge 0$ under either case:
\begin{equation}\label{bound_shd}
\sum_{n=1}^{\infty} \exp\left( \shd(\lambda_n,\sigma)\,t \right)
\le \exp\left( \shd(\lambda_1,\sigma)\,t\right) +\frac{\exp(c_{\mu_u}t)}{2\sqrt{\pi k t}}, 
\end{equation}
and 
\begin{equation}\label{bound_shl}
\sum_{n=1}^{\infty}  \shl(t, \lambda_n,\sigma)= \sum_{n=1}^{\infty} \int_0^t \exp\left( \shd(\lambda_n,\sigma)\,s \right) ds \leq  \frac{\exp\left( \shd(\lambda_1,\sigma) t\right)-1}{ \shd(\lambda_1,\sigma)} +\Theta(t, c_{\mu_u}),
\end{equation}
where 
\begin{equation}\label{bound_theta}
\Theta(t, c_{\mu_u}) = 
\begin{cases}
\displaystyle    \frac{1}{2\sqrt{k(-c_{\mu_u})}}, & c_{\mu_u}<0,\\[10pt]
\displaystyle  2 \sqrt{\frac{t}{\pi k}}\exp\left( c_{\mu_u} t\right) , &  c_{\mu_u} \ge 0.
\end{cases}
\end{equation}

In particular, for all $t \ge 0$,
\begin{equation}\label{bound_snh}
\sum_{n=1}^\infty \shn(t,\lambda_n,\sigma)
\le \shx(t,\lambda_1,\s)+\shy(t, c_{\mu_u}),
\end{equation}
with
\begin{align}
\shx(t,\lambda_1,\s) 
&= c_h \frac{\exp\left( \shd(\lambda_1,\sigma) t\right)-1}{ \shd(\lambda_1,\sigma)}
+ c_g\exp\left( \shd(\lambda_1,\sigma) t\right), \label{shx}\\[4pt]
\shy(t,c_{\mu_u})
&= c_h \Theta(t, c_{\mu_u})
+ c_g  \frac{\exp\left( c_{\mu_u}\,t\right)}{2\sqrt{\pi k t}}. \label{shy}
\end{align}
\end{lem}


\begin{proof}
We first prove \eqref{bound_shd} by splitting off $n=1$ and bound the tail.
The tail is bounded as follows.

Case A.  
Since $ \alpha_n^{u,\mu_u} \ge (n-1)\pi$, then for $n\ge2$,
\begin{align}\label{bnd_A}
& \sum_{n=2}^{\infty} 
\exp\left( \shd(\lambda_n,\sigma) t \right) = \sum_{n=2}^{\infty} 
\exp\left( (c_{\mu_u} - k\alpha_n^2) t \right)
\le 
\exp\left( c_{\mu_u}  t \right)
\sum_{m=1}^{\infty} 
\exp\left( -k\pi^2 m^2 t \right)\\
&\qquad \le 
\exp\left( c_{\mu_u} t \right)
\int_{0}^{\infty} 
\exp\left( -k\pi^2  x^2 t \right)\, dx
= 
\frac{\exp\left( c_{\mu_u} t \right)}{2\sqrt{\pi k t}}.
\end{align}

Case B. 
Since  $\alpha_n^{l,-\mu_u} \ge (n-\tfrac{1}{2})\pi$, then for $n\ge2$, 
\begin{align}\label{bnd_B}
&\sum_{n=2}^{\infty}  \exp\left( \shd(\lambda_n,\sigma) t \right) 
 \le 
\exp\left( c_{\mu_u} t \right)
\sum_{m=1}^{\infty} 
\exp\left( -k\pi^2 \left(m+\tfrac{1}{2}\right)^2 t \right)\\
&\qquad \le 
\exp\left( c_{\mu_u} t \right)
\int_{0}^{\infty} 
\exp\left( -k\pi^2  x^2 t \right)\, dx = 
\frac{\exp\left( c_{\mu_u} t \right)}{2\sqrt{\pi k t}}.
\end{align}

Finally, by \eqref{bnd_A} and \eqref{bnd_B}, the bound \eqref{bound_shd} is proved.

Now, we prove \eqref{bound_shl}.
Observe that
\begin{equation}\label{eq_521a}
\sum_{n=1}^{\infty} \shl(t,\lambda_n,\sigma)
=\sum_{n=1}^{\infty}\int_{0}^{t} \exp\left(\shd(\lambda_n,\sigma) s\right)ds
=\sum_{n=1}^{\infty} \frac{\exp\left( \shd(\lambda_n,\sigma) t\right)-1}{ \shd(\lambda_n,\sigma)}.
\end{equation}

\emph{Case 1: $c_{\mu_u}<0$.} We  split the case $n=1$ to the case $n> 1$ and start with the latter.
Since $\lambda_n$ are decreasing and $\shd(\lambda_n,\sigma) \leq c_{\mu_u}$, also $\shd(\lambda_n,\sigma) <0$ so, for $n\geq2$,
\begin{equation}\label{eq_521b}
\sum_{n=2}^{\infty} \frac{\exp\left( \shd(\lambda_n,\sigma) \ t\right)-1}{ \shd(\lambda_n,\sigma)} \leq 
\sum_{n=2}^{\infty} -\frac{1}{ \shd(\lambda_n,\sigma)}.
\end{equation}

Case A.  
For  $\alpha_n = \alpha_n^{u,\mu_u} $ the following holds
\begin{equation}\label{eq_525}
\sum_{n=2}^{\infty} -\frac{1}{\shd(\lambda_n,\sigma)} = \sum_{n=2}^{\infty} \frac{1}{k\alpha_n^2 - c_{\mu_u}}
\le \sum_{n=2}^{\infty} \frac{1}{k\pi^2(n-1)^2 - c_{\mu_u}}
= \sum_{m=1}^{\infty} \frac{1}{k\pi^2 m^2 - c_{\mu_u}} .
\end{equation}

Case B.
 For  $\alpha_n^{l,-\mu_u} $ the following holds
\begin{equation}\label{eq_525_b}
\sum_{n=2}^{\infty} -\frac{1}{\shd(\lambda_n,\sigma)} 
\le \sum_{n=2}^{\infty} \frac{1}{k\pi^2(n-\um)^2 - c_{\mu_u}}
= \sum_{m=1}^{\infty} \frac{1}{k\pi^2 (m+\um)^2 - c_{\mu_u}} <  \sum_{m=1}^{\infty} \frac{1}{k\pi^2 m^2 - c_{\mu_u}}.
\end{equation}

Thus, in either case, since $c_{\mu_u}\leq 0$, 
\begin{align}
&\sum_{m=1}^{\infty} \frac{1}{k\pi^2 m^2 - c_{\mu_u}}
\leq \frac{1}{k\pi^2} \sum_{m=1}^{\infty} \frac{1}{m^2 - \tfrac{c_{\mu_u}}{k\pi^2}} 
\leq \frac{1}{k\pi^2} \int_0^\infty \frac{dx}{x^2 - \tfrac{c_{\mu_u}}{k\pi^2}}\\
& =  \frac{1}{k\pi^2} \cdot \frac{1}{\sqrt{\tfrac{-c_{\mu_u}}{k\pi^2}}}
 \left[ \arctan \left(\frac{x}{\sqrt{\tfrac{-c_{\mu_u}}{k\pi^2}}}\right) \right]_{0}^{\infty}  =  \frac{1}{k\pi^2} \cdot \frac{1}{\sqrt{\tfrac{-c_{\mu_u}}{k\pi^2}}} \cdot \frac{\pi}{2} 
= \frac{1}{2\sqrt{k(-c_{\mu_u})}}.
\end{align}

After adding the first term of the series \eqref{eq_521a} to the above, the case $c_{\mu_u} <0$ of  \eqref{bound_shl} follows. 

\emph{Case 2: $ c_{\mu_u}  \geq 0.$}
Since the integrand in the second term of \eqref{eq_521a} is non negative, we can use Tonelli's theorem to invert the sum with the integral. 
We then split the case $n=1$ from the case $n> 1$. 
For the latter,
\begin{align}\label{eq_bnd5}
\sum_{n=2}^{\infty} \shl(t,\lambda_1,\s) = \int_{0}^{t} \exp  (c_{\mu_u}s) \left(\sum_{n=2}^{\infty}  \exp\left(  - k\alpha_n^2 s\right)\right)ds,
\end{align}
and as in \eqref{bnd_A} and \eqref{bnd_B}, use $\sum_{n=2}^{\infty}  \exp\left(  - k\alpha_n^2 s\right) \leq \int_{0}^{\infty} \exp(-k \pi^2 u^2 s )du <  \frac{1}{2\sqrt{k\pi s}}$ to obtain that 
\begin{equation}\label{eq_bnd6}
   \int_{0}^{t} \exp  (c_{\mu_u}s) \left(\sum_{n=2}^{\infty}  \exp\left(  - k\alpha_n^2 s\right)\right)ds \le \um \sqrt{\frac{1}{k\pi}}\int_{0}^{t}  \exp\left(c_{\mu_u} s\right)s^{-1/2} ds \le  \sqrt{\frac{1}{k\pi}} \exp\left(c_{\mu_u}  t\right)\sqrt{t},
\end{equation}
%
Adding the term $n=1$ of \eqref{eq_521a} to \eqref{eq_bnd6}
yields the case $c_{\mu_u}\ge 0$ of \eqref{bound_shl}.
Finally, \eqref{bound_snh} follows by \eqref{N_def_a}.
\end{proof}

A weaker bound that depends solely on the coefficients of the PDE, and not on $\lambda_1$, is provided by the following corollary.

\begin{cor}\label{bound_simple}Under the same assumptions of Lemma \ref{lem_bound_1} the following bounds hold:
\begin{equation}\label{bound_shd_1}
\sum_{n=1}^{\infty} \exp\left( \shd(\lambda_n,\sigma)\,t \right)
\le \exp\left(  c_{\mu_u}\,t\right)\left(1+\frac{1}{2 \sqrt{\pi k t}}\right), 
\end{equation}
and 
\begin{equation}\label{bound_shl_1}
\sum_{n=1}^{\infty}  \shl(t, \lambda_n,\sigma)\leq  \frac{\exp\left( c_{\mu_u} t\right)-1}{ c_{\mu_u}} + \Theta(t, c_{\mu_u}).
\end{equation}
\end{cor}
\begin{proof}
 The functions of $\shd(\lambda_1,\sigma)$ involved in the the first mode of \eqref{bound_shd} and  \eqref{bound_shl}  are strictly increasing: $e^{x}$, $\frac{e^{x}-1}{x}$. Since $c_{\mu_u} \geq \shd(\lambda_1,\sigma)$, this concludes the proof.
\end{proof}

Proposition \ref{prop_cn_bound} provides an estimate of  $c_{n,0}^{u,\mu}$ and $c_{n,\zeta}^{l,\mu}$ which is uniform in $n$,  and allows us to factor $c_{n,0}^{u,\mu_u}$ (resp. $c_{n,\zeta}^{l,-\mu_u}$) outside  the sum in  \eqref{dv_bound_1} (resp. \eqref{dv_bound_2}). Lemma \ref{lem_bound_1}, in turn, establishes a bound on the sums of $\shn$ in \eqref{dv_bound_1} and \eqref{dv_bound_2}. We can therefore combine these two results in the following theorem.

The first eigenvalue $\lambda_1^{u,\mu_u}$ has a structurally different form depending on  the sign of $\frac{\sigma^2}{\mu_u \zeta} $. We tackle the two cases separately.

\begin{thm}\label{316_bound_a}
Assume Hypothesis \ref{hyp} is verified. For any $t\geq 0$ and $x\in  \bar \shi$, the following holds.
\begin{equation}\label{316_eq}
\left| \frac{\partial V}{\partial x}(t,x) \right|  \leq  2 \frac{x}{\zeta} \frac{c^u_a }{ c^d_a } \exp\left(-\frac{\mu_u}{\sigma^2} x\right) \shq_0(t,\mu_u, \sigma),
\end{equation}
where

\begin{subequations}
\begin{equation}
\shq_0(t,\mu_u,\sigma)=
\end{equation}
\begin{empheq}[left={\empheqlbrace}]{align}
& c_{\alpha,0}^{u,\mu_u} \left(\shx(t,\lambda_1^{u,\mu_u},\s) +\shy(t,c_{\mu_u}) \right)
&& \mu_u \ge 0,\ \frac{\sigma^2}{\mu_u \zeta} \ge 1, \label{eq_380a}\\[6pt]
& \frac{\sigma^2 v_{\mu_u}}{\sigma^2\mu_u-(\mu_u^2-v_{\mu_u}^2)\zeta}\, \sinh \left(\frac{v_{\mu_u}}{\sigma^2}\zeta \right) \shx(t,\lambda_1^{u,\mu_u},\s) \\
&\qquad +\shy(t,c_{\mu_u})
&& \mu_u \ge 0,\ \frac{\sigma^2}{\mu_u \zeta} < 1, \label{eq_380}\\[6pt]
&\shx(t,\lambda_1^{l,-\mu_u},\s) +\shy(t,c_{-\mu_u}) 
&& \mu_u < 0, \label{eq_380aa_m}
\end{empheq}
\end{subequations}
\medskip
where 
  $c_{\alpha,0}^{u,\mu}$ is introduced in Proposition \ref{prop_cn_bound}, and $v_{\mu_u}$ in \eqref{v_mu_def} and  $\shx, \shy$ in \eqref{shx} and \eqref{shy}.

\end{thm}

\begin{proof}[Proof of Theorem \ref{316_bound_a}]

Let $x\in \shi$.

\textit{Case $\mu_u \geq 0 $ and $\frac{\sigma^2}{\mu_u \zeta} \geq 1$}. The conclusion follows after applying Item (1) of Proposition \ref{prop_cn_bound}, with $y=0$, to \eqref{dv_bound_1} and using \eqref{bound_snh}  to bound $\sum_{n=1}^\infty  \shn(t,\lambda^{u,\mu_u}_n,\sigma)$.

\textit{Case $\mu_u \geq 0 $ and $\frac{\sigma^2}{\mu_u \zeta} < 1$}. By applying Item (2)(a) of Proposition \ref{prop_cn_bound}, with $y = 0$ to \eqref{dv_bound_1} and considering only the terms with $n\ge2$, we obtain
\begin{equation}\label{eq352}
\sum_{n=2}^\infty \left| c^{u,\mu_u}_{n,0}(x) \right| \shn(t,\lambda^{u,\mu_u}_n,\sigma) \leq  2  \frac{ x}{\zeta} \exp\left(-\frac{\mu_u}{\sigma^2} x\right) \sum_{n=2}^\infty \shn(t,\lambda^{u,\mu_u}_n,\sigma).
\end{equation}
We then apply \eqref{bound_snh} observing that the sum is from $n =2$ hence it does not contain the first mode, that is $$\sum_{n=2}^\infty \shn(t,\lambda^{u,\mu_u}_n,\sigma) \leq \shy(t,c_{\mu_u}).$$
For $n=1$, we use Item (2)(b) of Proposition \ref{prop_cn_bound}, with $y = 0$, and \eqref{bound_snh}, considering only the first mode, 
\begin{align}
& \left| c^{u,\mu_u}_{1,0}(x) \right| \shn(t,\lambda^{u,\mu_u}_1,\sigma)\\
& \leq 2 \frac{x}{\zeta} \exp\left(-\frac{\mu_u}{\sigma^2} x\right) \sinh\left(\frac{v_{\mu_u}}{\s^2}\zeta\right) \frac{\s^2 v_{\mu_u}}{\s^2\mu_u -(\mu_u^2 -v_{\mu_u}^2) \zeta}  \shx(t,\lambda_1^{u,\mu_u},\s).\label{eq_363}
\end{align}
The conclusion follows after replacing \eqref{eq352} and \eqref{eq_363} inside \eqref{dv_bound_1}.

\textit{Case $\mu_u < 0 $}. Follows by using Item (3) of Proposition \ref{prop_cn_bound}, with $y = \zeta$, in \eqref{dv_bound_2} and \eqref{bound_snh}  to bound $\sum_{n=1}^\infty  \shn(t,\lambda^{l,-\mu_u}_n,\sigma)$.


For $x\in \partial \shi$, the conclusion follows by taking into account the boundary conditions \eqref{freeB2cc}, which we assume in Remark \ref{chg_var} are homogeneous.
\end{proof}


The following Proposition is similar to  Proposition \ref{prop_g_bound} but establishes, together with Theorem \ref{prop_d_bound}, estimates which are accurate near $x=\zeta$.

\begin{prop}\label{prop_d_bound}
Assume Hypothesis \ref{hyp} is in force. 
 For any $t\geq 0$ and $x\in  \shi$, the following holds.

\begin{enumerate}[(1)]
\item If $\mu_l \geq 0$, then
\begin{equation}\label{eq341_d}
\left| \frac{\partial V}{\partial x}(t,x) \right| \leq \frac{c^u_a }{ c^d_a }    \sum_{n=1}^\infty  \left| c^{l,\mu_l}_{n,\zeta}(x) \right| \shn(t,\lambda^{l,\mu_l}_n,\sigma), 
\end{equation}
\item otherwise, if $\mu_l < 0$:
\begin{equation}\label{eq341_d1}
\left| \frac{\partial V}{\partial x}(t,x) \right| \leq \frac{c^u_a}{ c^d_a }    \sum_{n=1}^\infty  \left| c^{u,-\mu_l}_{n,0}(\zeta-x) \right| \shn(t,\lambda^{u,-\mu_l}_n,\sigma), 
\end{equation}
\end{enumerate}
where $\lambda^{l,\mu_l}_n$ (resp. $\lambda^{u,-\mu_l}_n$) is as in \eqref{lambda_l_def} (resp. \eqref{lambda_u_def}) and $\shn(t,\lambda,\sigma)$ is as in \eqref{N_def_a}. Finally, $c^{l,\mu_l}_{n,\zeta}$ (resp. $c^{u,-\mu_l}_{n,0}$) are introduced in \eqref{cnDef2} (resp. \eqref{cnDef1}).
\end{prop}

\begin{proof}
We prove the result for $g=0$, the general case follows easily.

Item (1). By Lemma \ref{LThetaBound} and \eqref{bound_tauXd},
\begin{align}\label{eq342_d}
\left| \frac{\partial V}{\partial x}(t,x) \right| \leq  \frac{c^u_a \, c_{h}}{ c^d_a} \, \int_0^{t}  \exp\left(\sigma \, c_D \, c^u_a  \,   s\right) \P\left(\tau_\zeta^{\mu_l,x,S} \geq (c^d_a)^2 s \right) ds.
\end{align}
We now use \eqref{tauBound}, where $c_{n,\zeta}^{l,\mu_l}$ are defined in \eqref{cnDef2}, with $y = \zeta$, to compute the spectral expansion of $\P\left(\tau_\zeta^{\mu_l,x,S} \geq (c^d_a)^2 s \right)$.
By so doing, \eqref{eq342_d} becomes
\begin{align}\label{eq343_d}
\left| \frac{\partial V}{\partial x}(t,x) \right| &\leq  \frac{c^u_a \, c_{h}}{ c^d_a} \, \left| \int_0^{t} \left( \exp\left(\sigma \, c_D \, c^u_a  \,   s\right)  \sum_{n=1}^\infty  c_{n,\zeta}^{l,\mu_l}(x) \exp\left(-\lambda^{l,\mu_l}_n \, (c^d_a)^2 s\right) \right)ds\right|.
\end{align}
We now show that we can invert the sum with the integral.
For $n\geq 1$, by \eqref{lambda_l_def}, $\lambda^{l,\mu_l}_n = \frac{(\alpha^{l,\mu_l}_n)^2\s^2}{2 \zeta^2} + \frac{\mu_l^2}{2\s^2}$, where $\alpha^{l,\mu_l}_n \in ]\pim +(n-1)\pi, \pi +(n-1)\pi]$. This fact, together with Item (3) of Proposition \ref{prop_cn_bound}, yields
\begin{equation}\label{eq_517a}
\left| c_{n,\zeta}^{l,\mu_l}(x) \exp\left(-\lambda^{l,\mu_l}_n \, (c^d_a)^2 s\right) \right| \leq    \frac{ \zeta-x}{\zeta} \exp\left(-\frac{\mu}{\sigma^2} (x-\zeta)\right) \exp\left(-\frac{ (n-\um)^2\pi^2 \sigma^2 }{2\zeta^2} \, (c^d_a)^2 s\right).
\end{equation}
Therefore, by \eqref{eq_517a}, the Lebesgue dominated convergence theorem applies and \eqref{eq343} becomes
\begin{equation}\label{eq343aa}
\left| \frac{\partial V}{\partial x}(t,x) \right| \leq  \frac{c^u_a \, c_{h}}{ c^d_a} \, \left|  \sum_{n=1}^\infty \left( c_{n,\zeta}^{l,\mu_l}(x)  \int_0^{t}\exp\left((\sigma \, c_D \, c^u_a  -\lambda^{l,\mu_l}_n \, (c^d_a)^2 )   s\right)   ds \right)\right|.
\end{equation}
After resolving the integral, the proof is completed.

Item (2). By the reflection principle of Remark \ref{mu_neg_rem}, 
$$\P\left(\tau_\zeta^{\mu_l,x,S} \geq (c^d_a)^2 s \right)  = \sum_{n=1}^\infty  c_{n,0}^{u,-\mu_l}(\zeta-x) \exp\left(-\lambda^{u,-\mu_l}_n \, (c^d_a)^2 s\right), $$
then the result is proved as in Item (1). Use the same argument employed to derive \eqref{eq_510a} to apply the Lebesgue dominated convergence theorem. 
\end{proof}

Next result is similar to Theorem \ref{316_bound_a}. If $\mu_l \geq 0$, by applying the upper bound \eqref{bound_tauXd}, we do not need to consider two cases, depending on the values of $\frac{\s^2}{\mu_l \zeta}$.  This is because all the eigenvalues are of the form \eqref{lambda_def1}. 

\begin{rem}
The same conclusion of Lemma \ref{lem_bound_1} holds true if we replace $\mu_u$ with $\mu_l$ and  for $\lambda_n = \lambda_n^{u,-\mu_l}$ if $\mu_l < 0$ (Case  A) or  $\lambda_n = \lambda_n^{l,\mu_l}$ if $\mu_l \geq 0$ (Case B).
\end{rem}

\begin{thm}\label{317_bound}
Assume Hypothesis \ref{hyp} is in force. 
 For any $t\geq 0$ and $x\in \bar  \shi$, the following holds.

\begin{equation}
\left| \frac{\partial V}{\partial x}(t,x) \right|  \leq  2 \frac{\zeta-x}{\zeta} \frac{c^u_a}{ c^d_a } \exp\left(-\frac{\mu_l}{\sigma^2} (x-\zeta)\right) \shq_1(t,\mu_l, \sigma),
\end{equation}
where 
\begin{subequations}
\begin{equation}
\shq_1(t,\mu_l, \sigma)=
\end{equation}
\begin{empheq}[left={\empheqlbrace}]{align}
&\shx(t,\lambda_1^{l,\mu_l},\s)+\shy(t,c_{\mu_l})
 &&\mu_l \geq 0,  \label{eq_380aa}\\[6pt] 
& c_{\alpha,0}^{u,-\mu_l}  \left( \shx(t,\lambda_1^{u,-\mu_l},\s) +\shy(t,c_{-\mu_l})  \right) ,
&&\mu_l < 0 ,\ \frac{\sigma^2}{-\mu_l \zeta} \geq 1, \label{eq_380a_m}\\[6pt]
&  \frac{ \s^2 v_{-\mu_l}}{-\s^2\mu_l- (\mu_l^2 - v_{-\mu_l}^2) \zeta} \sinh \left(\frac{v_{-\mu_l}}{\s^2}\zeta\right)  \shx(t,\lambda_1^{u,-\mu_l},\s) \\
&\qquad + \shy(t,c_{-\mu_l})
&&\mu_l < 0 ,\ \frac{\sigma^2}{-\mu_l \zeta} < 1,\label{eq_380_m}
\end{empheq}
\end{subequations}
where 
  $c_{\alpha,0}^{l,-\mu_l}$ is introduced in Proposition \ref{prop_cn_bound}, and $v_{-\mu_l}$ in \eqref{v_mu_def} and  $\shx, \shy$ in \eqref{shx} and \eqref{shy}.

\end{thm}


\begin{proof}[Proof of Theorem \ref{317_bound}]
Let $x\in \shi$.

\textit{Case $\mu_l \geq 0$}. The conclusion follows after applying Item (3) of Proposition \ref{prop_cn_bound}, with $y=\zeta$, to \eqref{eq341_d} and using \eqref{bound_snh} to bound $ \sum_{n=1}^\infty  \shn(t,\lambda_n^{l, \mu_l},\sigma)$.

\textit{Case $\mu_l < 0$  and  $\frac{\sigma^2}{-\mu_l \zeta} \geq 1$}.   The conclusion follows after applying Item (1) of Proposition \ref{prop_cn_bound}, with $y=0$, to \eqref{eq341_d1} and  and using \eqref{bound_snh} to bound $\sum_{n=1}^\infty  \shn(t,\lambda^{u,-\mu_l}_n,\sigma) $.

\textit{Case $\mu_l < 0$  and  $\frac{\sigma^2}{-\mu_l \zeta} <1$}. 
For $n>1$, the result follows by applying Item (2)(a) of Proposition \ref{prop_cn_bound}, with $y = 0$,  to \eqref{eq341_d1}:
\begin{equation}\label{eq352_1}
\sum_{n=2}^\infty \left| c^{u,-\mu_l}_{n,0}(\zeta-x) \right| \shn(t,\lambda^{u,-\mu_l}_n,\sigma) \leq  2  \frac{\zeta- x}{\zeta} \exp\left(\frac{-\mu_l}{\sigma^2} (x-\zeta)\right) \sum_{n=2}^\infty \shn(t,\lambda^{u,-\mu_l}_n,\sigma).
\end{equation}
We then apply \eqref{bound_snh} observing that the sum in \eqref{eq352_1} is from $n =2$ hence it does not contain the first mode, that is $$\sum_{n=2}^\infty \shn(t,\lambda^{u,-\mu_l}_n,\sigma) \leq\shy(t,c_{-\mu_l}).$$
For $n=1$, by Item (2)(b) of Proposition \ref{prop_cn_bound}, with $y = 0$, and \eqref{bound_snh}, considering only the first mode, we obtain
\begin{align}
& \left| c^{u,-\mu_l}_{1,0}(\zeta-x) \right| \shn(t,\lambda^{u,-\mu_l}_1,\sigma)\\
& \leq 2 \frac{\zeta-x}{\zeta} \exp\left(-\frac{\mu_l}{\sigma^2} (x-\zeta)\right) \sinh\left(\frac{v_{-\mu_l}}{\s^2}\zeta\right) \frac{\s^2 v_{-\mu_l}}{-\s^2\mu_l -(\mu_l^2 -v_{-\mu_l}^2) \zeta} \,  \shx(t,\lambda_1^{u,-\mu_l},\s)   .\label{eq_363_1}
\end{align}
For $x\in \partial \shi$, the conclusion follows as in Theorem \ref{316_bound_a}.
\end{proof}

At this point, we are ready to prove the main result of the work. 

\begin{proof}[Proof of Theorem \ref{final_thm}]
For any $(t,x)$ both the bound of Theorem \ref{316_bound_a} and the bound of Theorem \ref{317_bound} hold.  
\end{proof}

As a straightforward corollary, we also provide a sufficient condition ensuring that   $\left| \frac{\partial V}{\partial x}(t,x) \right| $ remains bounded as $t \ra \infty$. This condition depends only on the first eigenvalue of the spectral decomposition \eqref{tauBound} and on  the coefficients of the problem. It is also very easy to verify in practice.

We remark that, in the case of drift $b(x)$ independent of time,  it is well known that $V$ converges,  as $t\ra \infty$,  to the solution of an ergodic PDE. This result holds even for fully non-linear parabolic PDE (see \cite{BDL}) and  in the semilinear case explicit convergence rates are provided in \cite{HM}. Consequently, in these settings, $ \lim_{t\ra \infty} \frac{\partial V}{\partial x}(t,x)$ exists.

When the drift depends on time, there are also results in this direction, provided that  $b(t,x) \ra b(x)$ as $t\ra \infty$, where $b(x)$ is a function independent of time. For example,  in  \cite[Chap. 11]{CC}  it is proved that, for $b(t,x)$ of the form $ \frac{b(x)}{t}$,  the solution $V(t,x)$ converges with rate $t^{-1}$ to the solution of an ergodic PDE.
 
In the general case, where $b(t,x)$ does not converge to any time-independent function, it is, to the best of our knowledge, not possible to establish the existence of a limit. Nevertheless, one can still prove that the solution remains bounded for all times. 

Recall $c_D$ is defined in \eqref{cost_def} and $c_a^u, c_a^d$ in Hypothesis \ref{hyp}. Also,  $\shd(\lambda,\sigma) = \sigma \, c_D \, c_a^u-\lambda  \, (c_a^d)^2$ was introduced in \eqref{D_def_a} and $\lambda^{l,\mu}_n$ (resp. $\lambda^{u,\mu}_n$) is as in \eqref{lambda_l_def} (resp. \eqref{lambda_u_def}).

\begin{cor}\label{lim_infty}
Let  $x\in  \bar \shi$. If the following are verified
\begin{align}
&\shd(\lambda_1^{u,\mu_u},\sigma)<0 \text{ or } \shd(\lambda_1^{l,\mu_l},\sigma)<0 \ \text{ if } \mu_u,\mu_l \geq 0, \label{der_infty_1}\\
&\shd(\lambda_1^{u,\mu_u},\sigma) <0 \text{ or } \shd(\lambda_1^{u,-\mu_l},\sigma) <0 \ \text{ if } \mu_u \geq 0,\mu_l < 0,\\
&\shd(\lambda_1^{l,-\mu_u},\sigma) <0 \text{ or } \shd(\lambda_1^{u,-\mu_l},\sigma) <0 \ \text{ if } \mu_u,\mu_l < 0,
\end{align}
then $\limsup_{t\ra \infty}\left| \frac{\partial V}{\partial x}(t,x) \right| <\infty$.
\end{cor}

\begin{rem}
If $c_{\mu_u}<0$ or $c_{\mu_l}<0$,  Corollary \ref{lim_infty} follows from Corollary \ref{bound_simple}. Such stronger hypothesis is easier to verify as it only depends on the coefficients of the problem and, in particular, does not depend on the first eigenvalue of the spectral expansion. 
\end{rem}

\begin{proof}
We prove the case $\mu_u,\mu_l \geq 0$, as the other items are verified in the same way.
By Theorem \ref{final_thm}, we shall prove that either $\lim_{t \ra \infty} \shq_0(t,\mu_u,\sigma) <\infty$ or $\lim_{t \ra \infty} \shq_1(t,\mu_l,\sigma) <\infty$.
This, in turn, according to \eqref{eq_380a}, \eqref{eq_380}, and \eqref{eq_380aa}  boils down to prove that $$\lim_{t\ra \infty} \sum_{n=1}^\infty \shn(t, \lambda^{u,\mu_u}_n,\sigma)<\infty, \quad \text{ or } \quad  \lim_{t\ra \infty} \sum_{n=1}^\infty \shn(t, \lambda^{l,\mu_l}_n,\sigma)<\infty.$$

We now show that $\lim_{t\ra \infty} \sum_{n=1}^\infty \shn(t, \lambda^{u,\mu_u}_n,\sigma)<\infty$, since the other case is similar. 
By assumption, $\lim_{t\ra \infty} \shn(t, \lambda^{u,\mu_u}_1,\sigma) = \frac{c_h}{-\shd(\lambda_1^{u,\mu_u},\s)} <\infty$.

For $n\geq 2$, since the eigenvalues are increasing, $\shl(t, \lambda^{u,\mu_u}_n,\sigma) \leq \frac{1}{-\shd(\lambda_1^{u,\mu_u},\s)}$. 
By \eqref{lambda_u_def} and recalling that $\alpha_n^{u,\mu_u}\in [(n-1)\pi, (n-1)\pi +\pim[$,  then   $\shd(\lambda_n^{u,\mu_u},\s) \leq f(n)$,  where $f(n) = \s c_D c_a^u -  \frac{(n-1)^2 \pi^2}{2\zeta^2} -\frac{\mu^2}{2\s^2}<0$. It follows that $\shl(t, \lambda^{u,\mu_u}_n,\sigma) \leq \frac{1}{-f(n)}  = O( \frac{1}{(n-1)^2})$.

Now, the other term of $\shn$ is $\exp\left(\shd(\lambda^{u,\mu_u}_n,\sigma)t\right)$. By assumption, for $n\geq 1$, it converges to $0$ as $t\ra \infty$ and, for $n\geq 2$, 
$$\sum_{n=2}^\infty \exp\left(\shd(\lambda^{u,\mu_u}_n,\sigma)t\right) \leq \sum_{n=2}^\infty \exp\left(f(n)t\right).$$ Moreover, since $f(n) = O( -(n-1)^2) $, Lebesgue dominated convergence applies and $$\lim_{t\ra \infty} \sum_{n=2}^\infty \exp\left(\shd(\lambda^{u,\mu_u}_n,\sigma)t\right) = \sum_{n=2}^\infty \lim_{t\ra \infty} \exp\left(\shd(\lambda^{u,\mu_u}_n,\sigma)t\right) = 0.$$
\end{proof}

%
%
%

\section{Appendix}

We prove two comparison theorems for reflected SDEs.
Since we were unable to find any results in the literature that handle time-dependent coefficients, we establish two versions which can be applied to equation \eqref{FC}.

\begin{prop}\label{comp}
Assume Hypothesis \ref{hyp}. Let $X$ be the solution to \eqref{FC}. If $x_1 \geq x_2$ then  a.s., for all $s\in [0,t]$,
\begin{align}
X_{s}^{t,x_1} \geq X_{s}^{t,x_2}.
\end{align}
\end{prop}

\begin{proof}
Let $\phi(x)= x^4 1_{\{x\geq 0\}}$. This function is $C^2$ and  non-negative  together with its first and second derivatives. In addition  
\begin{align}
|x|\phi'(x)&=x\phi'(x) = 4 \phi(x),\label{4.39}\\ 
x^2\phi''(x) &= 12 \phi(x). \label{4.40}
\end{align}
These facts will be useful later on. By It\^{o}'s formula,
\begin{align}
&\phi(X_s^{t,x_2}- X_s^{t,x_1})= \phi(x_2-x_1) \label{aa} \\
&\qquad + \int_0^s \phi' (X_r^{t,x_2}-X_r^{t,x_1}) \left(  b(t-r,X_r^{t,x_2}) -b(t-r,X_r^{t,x_1}) \right)dr  \label{aaa}\\
&\qquad + \int_0^s  \sigma \phi' (X_r^{t,x_2} - X_r^{t,x_1}) \left( a(X_r^{t,x_2})- a(X_r^{t,x_1}) \right) dB_r \label{bb} \\
&\qquad + \int_0^s \frac{\sigma^2}{2} \phi'' (X_r^{t,x_2}- X_r^{t,x_1}) \left( a(X_r^{t,x_2}) - a(X_r^{t,x_1}) \right)^2 dr \label{cc}\\
&\qquad + \int_0^s  \phi' (X_r^{t,x_2}- X_r^{t,x_1}) (dL^{t,x_2}_r -dL^{t,x_1}_r) - \int_0^s  \phi' (X_r^{t,x_2}- X_r^{t,x_1}) (dU^{t,x_2}_r -dU^{t,x_1}_r).
\end{align}
Expanding the last  line we obtain
\begin{align}
& \int_0^s \phi' (X^{t,x_2}_r- X_r^{t,x_1}) dL^{t,x_2}_r - \phi'(X^{t,x_2}_r -X^{t,x_1}_r)dL^{t,x_1}_r \\
& \qquad  - \left( \phi'(X^{t,x_2}_r-X^{t,x_1}_r) dU^{t,x_2}_r - \phi'(X^{t,x_2}_r -X^{t,x_1}_r)dU^{t,x_1}_r \right).\label{dd}
\end{align}
All of the four terms in the last integral are non positive: in fact when $dL^{t,x_2}_r$ is positive, then $X^{t,x_2}_r= 0$ so  $X^{t,x_1}_r \geq X^{t,x_2}_r $, and $\phi' (X^{t,x_2}_r -X^{t,x_1}_r) dL^{t,x_2}_r=0$. When $dU^{t,x_1}_r$ is positive, then  $X^{t,x_1}_r= \zeta$ and $ \phi' (X^{t,x_2}_r -X^{t,x_1}_r) dU^{t,x_1}_r = 0$. The remaining two terms are negative, hence we can get a larger value by replacing these terms  by 0. Using the Lipschitz property of $x\ra b(\cdot,x)$, we can find an upper bound to the term in  \eqref{aaa}, 
\begin{align}
&  \int_0^s \phi' (X_r^{t,x_2}-X_r^{t,x_1}) \left( b(t-r,X_r^{t,x_2}) - b(t-r,X_r^{t,x_1}) \right)dr  \label{dd0} \\
& \leq \int_0^s \phi' (X_r^{t,x_2}-X_r^{t,x_1})   C  \left| X_r^{t,x_2}-X_r^{t,x_1} \right| dr  \label{dd1}\\
&\leq  \int_0^s  C \phi (X_r^{t,x_2}-X_r^{t,x_1})  dr, \label{dd2}
\end{align}
where the last inequality \eqref{dd2} follows from the property \eqref{4.39} of $\phi$. We can therefore rewrite \eqref{aa}--\eqref{cc}, also using the Lipschitz property of $a$ and \eqref{4.40}, as
\begin{align}
&\phi(X_s^{t,x_2}-X_s^{t,x_1}) \leq   \phi(x_2-x_1)+ \int_0^s  C \phi (X_r^{t,x_2}-X_r^{t,x_1})  dr \\
&  \qquad + \int_0^s \sigma  \phi' (X_r^{t,x_2}-X_r^{t,x_1}) \left( a(X_r^{t,x_2})-  a(X_r^{t,x_1}) \right) dB_r .
\end{align}
Note that $  \phi(x_2-x_1)=0$. 
The stochastic integral is a martingale, hence taking expectations, the above becomes
\begin{align}
\E  \left( \phi (X_s^{t,x_2}-X_s^{t,x_1})\right) \leq C \E \left( \int_0^s  \phi (X_r^{t,x_2}-X_r^{t,x_1}) dr\right),
\end{align}
so, by Gronwall's Lemma (see  for example \cite[Chapter VII, Section X]{walter1998ode}), $$\E \left(\phi(X_s^{t,x_2}-X_s^{t,x_1})\right) =0.$$ Since $X$ is continuous, then  a.s. for all $s\in [0,t]$, $\phi(X_s^{t,x_2}-X_s^{t,x_1}) = 0$, i.e. $X_s^{t,x_2} \leq X_s^{t,x_1}$.
\end{proof}

The proposition below is another form of comparison theorem on the drift of the SDE.

\begin{prop}\label{comp1}
Assume Hypothesis \ref{hyp}.  Let $X^{t,x,1}$ be an SDE driven by $b_1 \geq b$, $b_1$ satisfying  Item (3) of Hypothesis \ref{hyp}. Then  a.s., for all $s\in [0,t]$, and $x \in \bar \shi$,
\begin{align}
X_{s}^{t,x,1} \geq X_{s}^{t,x}.
\end{align}
\end{prop}

\begin{proof}
We apply It\^{o}'s formula to  $\phi(X_s^{t,x}- X_s^{t,x,1})$. The proof 
uses the same arguments of the one of Proposition \ref{comp} after replacing \eqref{aaa} by
\begin{equation} \label{aaa1}
 \int_0^s \phi' (X_r^{t,x}-X_r^{t,x,1}) \left(  b(t-r,X_r^{t,x}) -b_1(t-r,X_r^{t,x,1}) \right)dr.
\end{equation}
Now, use the following upper bound, in place of \eqref{dd0}--\eqref{dd2}, to \eqref{aaa1} 
\begin{align}
&  \int_0^s \phi' (X_r^{t,x}-X_r^{t,x,1}) \left( b(t-r,X_r^{t,x}) - b_1(t-r,X_r^{t,x,1}) \right)dr\\
& = \int_0^s \phi' (X_r^{t,x}-X_r^{t,x,1}) \left( b(t-r,X_r^{t,x}) - b(t-r,X_r^{t,x,1}) \right)dr\\
&\qquad + \int_0^s \phi' (X_r^{t,x}-X_r^{t,x,1}) \left( b(t-r,X_r^{t,x,1}) - b_1(t-r,X_r^{t,x,1}) \right)dr  \label{cc11}\\
&\leq  \int_0^s  C \phi (X_r^{t,x_2}-X_r^{t,x_1})  dr, \label{dd22}
\end{align}
where in the last inequality we used the Lipschitz property of $x\ra b(\cdot,x)$ and  the fact that \eqref{cc11} is non-positive.
\end{proof}

\section*{Acknowledgements}

The author gratefully acknowledges the valuable discussions and suggestions provided by Professor R.C. Dalang during the development of this work.

\bibliographystyle{plain}
\bibliography{biblioS}

@preamble{"\def\cprime{$'$} "}

@article{n,
	Author = {Nourdin, N.},
	Journal = {Preprint},
	Title = {Titre \`a preciser}}

@book{ry3,
	Address = {Berlin},
	Author = {Revuz, D. and Yor, M.},
	Edition = {Third},
	Publisher = {Springer-Verlag},
	Series = {Grundlehren der Mathematischen Wissenschaften [Fundamental Principles of Mathematical Sciences]},
	Title = {Continuous martingales and {B}rownian motion},
	Volume = {293},
	Year = {1999}}

@article{z,
	Author = {Zvonkin, A. K.},
	Journal = {Mat. Sb. (N.S.)},
	Pages = {129--149, 152},
	Title = {A transformation of the phase space of a diffusion process that will remove the drift},
	Volume = {93(135)},
	Year = {1974}}

@book{Lunardi,
author="Lunardi, A.",
title="{Analytic semigroups and optimal regularity in parabolic problems.
    Reprint of the 1995 hardback ed.}",
language="English",
publisher="{Modern Birkh\"auser Classics. Basel: Birkh\"auser. xvii, 424~p.}",
year="2013",
classmath="{*35-02 (Research monographs (partial differential equations))

}",
}

@book{Fri,
author = "Friedman, A.",
title= {Partial Differential Equations of Parabolic Type},
language="English",
publisher={Dover},
year="2008. Reprint of the 1964 edition published by Prentice-Hall",
}

@article{DeuZa,
title = "Bismuth-{E}lworthy's formula and random walk representation for SDEs with reflection ",
journal = "Stochastic Processes and their Applications",
volume = "115",
pages = "907 -925",
year = "2005",
author = "Deuschel, J-D. and Zambotti, L.",
}

@article{Lin,
title = "Lookback options and diffusion hitting times: a spectral expansion approach ",
journal = "Finance and Stochastics",
volume = "8",
pages = "373-398",
year = "2004",
author = "Linetsky, V.",
}

@article{And,
title = "Pathwise differentiability for SDEs in a smooth domain with reflection",
journal = "Electronic Journal of Probability",
volume = "16",
pages = "845-879",
year = "2011",
author = "Sebastian, A.",
}

@book{LUS,
   title =     {Linear and Quasilinear Equations of Parabolic Type },
   author =    {Ladyzenskaja, O.A. and Solonnikov, V.A. and Uralceva, N.N.},
   publisher = {AMD},
   year =      {1968},
   edition =   {1}
}

@article{tanaka1979,
author = "Tanaka, H.",
fjournal = "Hiroshima Mathematical Journal",
journal = "Hiroshima Math. J.",
number = "1",
pages = "163--177",
publisher = "Hiroshima University, Department of Mathematics",
title = "Stochastic differential equations with reflecting boundary condition in convex regions",
volume = "9",
year = "1979"
}

@article{kim1990,
author = "Kim, J. and Pollard, D.",
doi = "10.1214/aos/1176347498",
fjournal = "The Annals of Statistics",
journal = "Ann. Statist.",
month = "03",
number = "1",
pages = "191--219",
publisher = "The Institute of Mathematical Statistics",
title = "Cube Root Asymptotics",
volume = "18",
year = "1990"
}

@book{CC,
  title={Optimal Solution and Asymptotic Properties of a Stochastic Control Problem Arising in Sailboat Trajectory Optimization},
  author={Ciccarella, C.},
  year={2017},
  publisher={EPFL PhD thesis},
}

@article{CDV1,
      title={An Optimal Control Problem Arising in Sailboat Trajectory Optimization}, 
      author={Ciccarella, C. and Dalang, R.C. and Vinckenbosch, L.},
      year={2024},
      eprint={2404.03773},
      archivePrefix={arXiv},
      primaryClass={math.OC},
      note={https://arxiv.org/abs/2404.03773}, 
}

@article{BMhit,
        title={The Hitting Time Density for a Reflected Brownian Motion},
        author={Hu, Q. and Wang, Y. and Yang, X.},
        journal={Computational Economics},
        publisher={Springer},
        year={2012},
        volume = {40},
        pages ={1-18}        
}

@article{Barles2021,
  author    = {Barles, G.},
  title     = {Local Gradient Estimates for Second-Order Nonlinear Elliptic and Parabolic Equations by the Weak {B}ernsteins Method},
  journal   = {Partial Differential Equations and Applications},
  year      = {2021},
  volume    = {2},
  number    = {6},
  pages     = {71}
}

@article{BartierSouplet2004,
  author    = {Bartier, J-P. and Souplet, P.},
  title     = {Gradient bounds for solutions of semilinear parabolic equations without Bernstein's quadratic condition},
  journal   = {Comptes Rendus Mathematique},
  year      = {2004},
  volume    = {338},
  number    = {7},
  pages     = {533--538},
  doi       = {10.1016/j.crma.2003.12.030},
  url       = {https://www.numdam.org/item/10.1016/j.crma.2003.12.030.pdf}
}

@article{Engler1990,
  author    = {Engler, H. and Kawohl, B. and Luckhaus, S.},
  title     = {Gradient Estimates for Solutions of Parabolic Equations and Systems},
  journal   = {Journal of Mathematical Analysis and Applications},
  year      = {1990},
  volume    = {147},
  number    = {2},
  pages     = {309--329},
  doi       = {10.1016/0022-247X(90)90340-8}
}

@article{LeBalch2021,
  author    = {Le Balc'h, K.},
  title     = {Exponential bounds for gradient of solutions to linear elliptic and parabolic equations},
  journal   = {Journal of Functional Analysis},
  volume    = {281},
  number    = {8},
  pages     = {109139},
  year      = {2021},
  doi       = {10.1016/j.jfa.2021.109139},
  url       = {https://www.sciencedirect.com/science/article/abs/pii/S0022123621001762}
}

@article{Kresin2020,
  author    = {Kresin, G. and Maz'ya, V.},
  title     = {Sharp pointwise estimates for the gradients of solutions to linear parabolic second-order equations in the layer},
  journal   = {Applicable Analysis},
  volume    = {99},
  number    = {10},
  pages     = {1767--1785},
  year      = {2020},
  doi       = {10.1080/00036811.2020.1732356},
  url       = {https://www.tandfonline.com/doi/abs/10.1080/00036811.2020.1732356}
}

@article{Han1998,
  author  = {Han, Q.},
  title   = {On the {S}chauder estimates of solutions to parabolic equations},
  journal = {Annali della Scuola Normale Superiore di Pisa - Classe di Scienze},
  volume  = {27},
  number  = {1},
  pages   = {1--26},
  year    = {1998},
  url     = {https://eudml.org/doc/84352}
}

@article{Kukavica2024,
title = {On the pointwise {S}chauder estimates for elliptic and parabolic equations},
journal = {Discrete and Continuous Dynamical Systems},
volume = {44},
number = {12},
pages = {3734-3759},
year = {2024},
issn = {1078-0947},
doi = {10.3934/dcds.2024076},
url = {https://www.aimsciences.org/article/id/666012fdf4772e1821a779b3},
author = {Kukavica, I. and Le, Q.}
}

@article{solonnikov2016,
  author    = {Solonnikov, V. A.},
  title     = {Proof of {S}chauder estimates for parabolic initial-boundary value model problems via O.A. {L}adyzhenskaya’s {F}ourier multipliers theorem},
  journal   = {Zap. Nauchn. Sem. POMI},
  volume    = {444},
  pages     = {133--156},
  year      = {2016},
  publisher = {Steklov Mathematical Institute of RAS},
}

@article{Montoro2018,
  author    = {Montoro, L. and Punzo, F. and Sciunzi, B.},
  title     = {Pointwise estimates for solutions of semilinear parabolic inequalities with a potential},
  journal   = {Rendiconti Lincei - Matematica e Applicazioni},
  volume    = {29},
  number    = {2},
  pages     = {255--288},
  year      = {2018},
  doi       = {10.4171/RLM/804},
  url       = {https://ems.press/journals/rlm/articles/15417},
}

@book{walter1998ode,
  author    = {Walter, W.},
  title     = {Ordinary Differential Equations},
  publisher = {Springer},
  year      = {1998},
  edition   = {Graduate Texts in Mathematics, Vol. 182},
  isbn      = {978-0387985134}
}

@misc{crisan2023,
  title={A Probabilistic Representation of the Derivative of a One-Dimensional Killed Diffusion Semigroup},
  author={Crisan, D. and Kohatsu-Higa, A.},
  year={2023},
  eprint={2312.07084},
  archivePrefix={arXiv},
  primaryClass={math.PR},
  note={https://arxiv.org/abs/2312.07084}
}

@book{Oksendal2003,
  author    = {Øksendal, B.},
  title     = {Stochastic Differential Equations: An Introduction with Applications},
  edition   = {6},
  publisher = {Springer},
  series    = {Universitext},
  year      = {2003},
  isbn      = {978-3-540-04758-0},
  doi       = {10.1007/978-3-642-14394-6},
  url       = {https://doi.org/10.1007/978-3-642-14394-6}
}

@article{CassFriz2006,
  author    = {Cass, T.R. and Friz, P.K.},
  title     = {The {B}ismut--{E}lworthy--{L}i formula for jump-diffusions and applications to {M}onte {C}arlo pricing in finance},
  year      = {2006},
  note       = {https://arxiv.org/abs/math/0604311}
}

@book{Glasserman1991,
  author    = {Glasserman, P.},
  title     = {Gradient Estimation via Perturbation Analysis},
  publisher = {Kluwer Academic Publishers},
  address   = {Boston, MA},
  year      = {1991},
  series    = {The Springer International Series in Engineering and Computer Science},
  volume    = {116},
  isbn      = {978-0-7923-9095-4}
}

@book{fleming2006,
  title     = {Controlled Markov Processes and Viscosity Solutions},
  author    = {Fleming, W. H. and Soner, H. M.},
  year      = {2006},
  edition   = {2},
  publisher = {Springer},
  address   = {New York},
  series    = {Stochastic Modelling and Applied Probability},
  volume    = {25},
  isbn      = {978-0-387-26045-7},
  doi       = {10.1007/0-387-31071-1},
  url       = {https://link.springer.com/book/10.1007/0-387-31071-1}
}

@article{BDL,
  author    = {Barles, G. and Da Lio, F.},
  title     = {On the boundary ergodic problem for fully nonlinear equations in bounded domains with general nonlinear Neumann boundary conditions},
  journal   = {Annales de l'Institut Henri Poincaré, Analyse non linéaire},
  volume    = {22},
  number    = {5},
  pages     = {521--541},
  year      = {2005},
  publisher = {Elsevier}
}

@article{HM,
  author    = {Hu, Y.  and Madec, P-Y},
  title     = {A probabilistic approach to large time behaviour of viscosity solutions of parabolic equations with Neumann boundary conditions},
  journal   = {Applied Mathematics and Optimization},
  volume    = {74},
  pages     = {345--374},
  year      = {2015}
}

\end{document}